\documentclass[centertags,leqno]{article}

\usepackage{amsmath}
\usepackage{amssymb}
\usepackage{latexsym}
\usepackage{amsxtra}
\usepackage{amscd}
\usepackage{theorem}

\usepackage{curves}
\usepackage{epic}
\usepackage{eepic}

\theoremstyle{plain}

{\theorembodyfont{\slshape}

        \newtheorem{thm}{Theorem}[section]
        \newtheorem{cor}[thm]{Corollary}
        \newtheorem{lem}[thm]{Lemma}
        \newtheorem{prop}[thm]{Proposition}
        
}

{\theorembodyfont{\rmfamily}

        \newtheorem{defn}[thm]{Definition}
        
        \newtheorem{rem}[thm]{Remark}

}

\renewcommand{\em}{\sl}

\newcommand{\proof}{{\bf Proof:\ }}
\newcommand{\Endproof}{\hspace*{\fill} $\Box$ \vspace{1ex} \noindent }

\makeatletter
\renewcommand{\subsection}{\@startsection{subsection}{2}%
        {\z@}{-3.25ex plus -1ex minus-.2ex}{-1em}{\bf}}
\makeatother

\newcommand{\PP}{\mathbb{P}}
\newcommand{\ZZ}{\mathbb{Z}}
\newcommand{\CC}{\mathbb{C}}
\newcommand{\RR}{\mathbb{R}}
\newcommand{\QQ}{\mathbb{Q}}

\newcommand{\FF}{\mathbb{F}}
\renewcommand{\AA}{\mathbb{A}}

\newcommand{\F}{\mathcal{F}}
\newcommand{\E}{\mathcal{E}}
\newcommand{\V}{\mathcal{V}}
\newcommand{\W}{\mathcal{W}}
\newcommand{\U}{\mathcal{U}}

\newcommand{\LL}{\mathcal{L}}
\newcommand{\OO}{\mathcal{O}}

\newcommand{\bC}{{\bf C}}

\newcommand{\g}{{\bf g}}
\newcommand{\p}{\mathfrak{p}}
\newcommand{\GL}{{\rm GL}}
\newcommand{\SL}{{\rm SL}}
\newcommand{\PGL}{{\rm PGL}}
\newcommand{\PSL}{{\rm PSL}}
\newcommand{\Gal}{{\rm Gal}}
\newcommand{\Aut}{{\rm Aut}}
\newcommand{\Spec}{{\rm Spec\,}}

\newcommand{\Ni}{{\rm Ni}}

\newcommand{\betat}{\tilde{\beta}}
\newcommand{\betab}{\bar{\beta}}
\newcommand{\etab}{\bar{\eta}}
\newcommand{\pib}{\bar{\pi}}
\newcommand{\varphib}{\bar{\varphi}}

\newcommand{\kb}{\bar{k}}
\newcommand{\Bb}{\bar{B}}

\newcommand{\Qb}{\bar{\QQ}}

\newcommand{\para}{_p}

\newcommand{\et}{{\rm\acute{e}t}}
\newcommand{\geo}{^{\rm geo}}
\newcommand{\tp}{^{\rm top}}
\newcommand{\an}{^{\rm an }}
\newcommand{\red}{^{\rm red }}

\newcommand{\sing}{_{\rm\scriptscriptstyle sing}}
\newcommand{\dR}{_{\rm dR}} 

\newcommand{\inj}{\hookrightarrow}
\newcommand{\To}{\;\longrightarrow\;}
\newcommand{\iso}{\stackrel{\sim}{\to}}
\newcommand{\liso}{\;\stackrel{\sim}{\longrightarrow}\;}

\newcommand{\gen}[1]{\mathopen\langle#1\mathclose\rangle}

\title{Variation of local systems and parabolic cohomology}

\author{
   Michael Dettweiler \\ Universit\"at Heidelberg
   \and
   Stefan Wewers
 \\ Universit\"at Bonn}

\date{}

\begin{document}

\maketitle

\begin{abstract}
  Given a family of local systems on a punctured Riemann sphere,
  with moving singularities, its first parabolic cohomology is a local
  system on the base space. We study this situation from different
  points of view. For instance, we derive universal formulas for the
  monodromy of the resulting local system. We use a particular example
  of our construction to prove that the simple groups $\PSL_2(p^2)$
  admit regular realizations over the field $\QQ(t)$ for primes
  $p\not\equiv 1,4,16\mod{21}$. Finally, we compute the monodromy of
  the Euler-Picard equation, reproving a classical result of Picard.
\end{abstract}

%------------------------------------------------------

\section*{Introduction}

%\subsection*{Variations of local systems and the middle convolution}

Local systems on the punctured Riemann sphere arise in various
branches of mathematics and have been intensively studied, see e.g.\ 
\cite{DeligneED}, \cite{KatzRLS}. One way to produce such local
systems with interesting properties is the following. Suppose we are
given a family $\V_s$ of local systems on the punctured sphere with
moving singularities, parameterized by some base space $S$. More
precisely, let $D\subset\PP^1\times S$ be a smooth relative divisor
and let $\V$ be a local system on $U:=\PP^1\times S-D$. Let
$\pib:\PP^1\times S\to S$ denote the second projection and
$j:U\inj\PP^1\times S$ the canonical inclusion. Then the first higher
direct image sheaf
\[
     \W \;:=\; R^1\pib_*(j_*\V)
\]
is a local system on $S$, whose stalk at a point $s\in S$ is the first
parabolic cohomology group of the local system $\V_s$, the restriction
of $\V$ to the fiber $U_s:=\pi^{-1}(s)$. Choose a base point $s_0\in
S$ and set $\V_0:=\V_{s_0}$. We call $\V$ a {\em variation} of the
local system $\V_0$ over the base $S$ and $\W$ the first parabolic
cohomology of the variation $\V$.

\vspace{2ex} A special case of this construction is the {\em middle
  convolution} studied by N.\ Katz in \cite{KatzRLS} and by Dettweiler
and Reiter in \cite{DettwReiterMC}. Starting with some local system
$\V_0$ on the punctured Riemann sphere $U_0$, one constructs a
variation of local systems over the base $U_0$ by twisting $\V_0$ with a
$1$-dimensional system with two singularities, one of which is moving
over all points of $U_0$. The parabolic cohomology of this variation
gives rise to a local system on $U_0$, called the {\em middle
  convolution} of $\V_0$. In \cite{KatzRLS}, Katz proves that all
rigid local systems can be constructed from $1$-dimensional systems by
successive application of middle convolution and `scaling'.

A local system on the punctured Riemann sphere corresponds to an
$r$-tuple of invertible matrices $\g=(g_1,\ldots,g_r)\in\GL_n(\CC)^r$.
V\"olklein \cite{Voelklein01} and, independently, Dettweiler and
Reiter \cite{DettwReiterKatz} have defined an operation
$\g\mapsto\tilde{\g}$ on tuples of invertible matrices over any field
$K$ corresponding to the middle convolution (in \cite{Voelklein01}, it
is called the {\em braid companion functor}).  The definition of this
operation needs only simple linear algebra.  Therefore, the tuple
$\tilde{\g}$ can be easily computed, whereas in the original work of
Katz the matrices $\tilde{g}_i$ are computed only up to conjugation in
$\GL_m$. This construction has had many applications to the Regular
Inverse Galois Problem, see e.g.\ \cite{Voelklein01} and
\cite{DettwReiterKatz}

\vspace{2ex} The goal of the present paper is to study the parabolic
cohomology of an arbitrary variation of local systems (see the
beginning of this introduction), both from an analytic and from an
arithmetic point of view.

In the first part of our paper, we treat the analytic aspect, i.e.\ we
use singular cohomology. Given a variation $\V$ of local systems over
a base $S$, we present an effective method to compute the monodromy
representation of $\pi_1(S)$ on the parabolic cohomology of $\V$. The
result depends on the tuple of matrices corresponding to the fibers of
the variation $\V$, and on the map from $\pi_1(S)$ to the Hurwitz
braid group which describes how the singularities of these fibres move
around on $\PP^1$. Formally, our method is a straightforward extension
of V\"olklein's braid companion functor. Apart from the greater
generality, the main difference to
V\"olklein's approach is that we provide a cohomological
interpretation of our computation. Such an interpretation has the
advantage that it makes it very easy, using comparison theorems, to
translate results from one world into another. For instance, one can
use topological methods to compute the monodromy of a local system.
Then one can apply the results of this computation either to the
Galois representations or to the differential equation attached to the
local system in question.

We give two applications which illustrate this principle. The first example is
concerned with the Regular Inverse Galois Problem and generalizes the method
of \cite{DettwReiterKatz} and \cite{Voelklein01}. Using this generalization,
we prove the following theorem.

\vspace{2ex}
\noindent{\bf Theorem:} {\em The simple group $\PSL_2(p^2)$ admits a regular
realization as Galois group over $\QQ(t)$, for all primes $p\not\equiv
1,4,16\mod{21}$.}

\vspace{2ex} The only other cases where the group $\PSL_2(p^2)$ is
known to admit a regular realization over $\QQ(t)$ are for
$p\not\equiv\pm 1\mod{5}$, by a result of Feit \cite{Feit85}, and for
$p\not\equiv\pm 1\mod{24}$, by a result of Shiina \cite{Shiina1},
\cite{Shiina2} (see also \cite{Mestre91}, \cite{DettweilerTokyo} and
\cite{VoelkleinCrelle}).

If $q<n$ then one knows that the group ${\rm PSp}_{2n}(q)$ occurs regularly
over $\QQ(t)$, see \cite{VoelkleinThompson} and \cite{DettwReiterKatz}.
Similar bounds exist also for other classical groups. On the other hand,
experience shows that it is much harder to realize classical groups of small
Lie rank. The realizations of $\PSL_2(p^2)$ in the above theorem all come from
one particular variation of local systems. It is very likely that, by choosing
different variations, one can realize many more series of classical groups of
small rank. We have chosen one particular example leading to the above
theorem, because the case of rank one seems to us the hardest case.

\vspace{2ex} In the last section, we compute the monodromy of the
Picard--Euler equation, reproving a classical result of Picard
\cite{Picard83},\cite{Picard84}. The Picard--Euler equation is the
Fuchsian system of partial differential equations associated to the
universal family of Picard curves (tricyclic covers of the Riemann
sphere with five branch points). One can identify the local system of
solutions to this equation with the parabolic cohomology of a
variation of local systems. Therefore, the method developed in the
first part of this paper can be used to compute the monodromy of the
Picard--Euler equation. We expect that the same approach will prove
useful in the study of more general differential equations. There is a
vast literature dealing with special classes of Fuchsian systems, many
of which arise as the parabolic cohomology of a variation. For
instance, this is the case for hypergeometric systems, see
\cite{KatzRLS}, \cite{DeligneMostow}. So far, there seemed to be no
general method available to compute the monodromy of a local system
coming from a variation. Here we present such a method. It is very
general, explicit and can easily be implemented on a computer.

\vspace{2ex} The first author would like to thank the mathematical
department of the University of Tel Aviv for their hospitality during
his stay in spring 2003, especially M. Jarden and D. Haran.  Both
authors acknowledge the financial support provided through the
European Community's Human Potential Program under contract
HPRN-CT-2000-00114, GTEM.

%------------------------------------------------------------------

\section{Parabolic cohomology} 
  \label{locals}

We study the first parabolic cohomology of a local system on the punctured
sphere. In particular, we show that it is isomorphic to a certain module
$W_\g$, defined in \cite{Voelklein01}. 

\subsection{} \label{locals0}

Let $X$ be a connected and locally contractible topological space.
Let $R$ be a commutative ring with unit.  A {\em local system} of
$R$-modules on $X$ is a locally constant sheaf $\V$ on $X$ whose
stalks are free $R$-modules of finite rank. We denote by $\V_x$ the
stalk of $\V$ at a point $x\in X$. If $f:Y\to X$ is a continuous map
then $\V_f$ denotes the group of global sections of the sheaf
$f^*\V$. Note that if $Y$ is simply connected then the natural
morphism
\[
    \V_f \;\To\; \V_{f(y)}
\]
is an isomorphism, for all $y\in Y$. Therefore, a path
$\alpha:[0,1]\to X$ gives rise to an isomorphism
\[
   \V_{\alpha(0)} \;\liso\; \V_{\alpha(1)},
\]
obtained as the composition of the isomorphisms $\V_{\alpha(0)}\cong\V_\alpha$
and $\V_\alpha\cong\V_{\alpha(1)}$. The image of $v\in\V_{\alpha(0)}$ under
the above isomorphism is denoted by $v^\alpha$. It only depends on the
homotopy class of $\alpha$. 

Let us fix a base point $x_0\in X$ and set $V:=\V_{x_0}$. We let elements of
$\GL(V)$ act on $V$ from the {\em right}. Then the map 
\[
       \rho:\pi_1(X,x_0) \;\To\; \GL(V),
\]
defined by $v\cdot\rho(\alpha):=v^\alpha$, is a group homomorphism,
i.e.\ a representation of $\pi_1(X,x_0)$. It is a standard fact that
the functor $\V\mapsto V:=\V_{x_0}$ is an equivalence of categories
between local systems on $X$ and representations of $\pi_1(X,x_0)$.

\subsection{} \label{locals1} 

Let $X$ be a compact (topological) surface of genus $0$ and $D\subset
X$ a subset of cardinality $r$. We set $U:=X-D$. There exists a
homeomorphism $\kappa:X\iso\PP^1_\CC$ between $X$ and the Riemann
sphere which maps the set $D$ to the real line $\RR\subset\PP^1_\CC$.
Such a homeomorphism is called a {\em marking} of $(X,D)$. 

Let us, for the moment, identify $X$ with $\PP^1_\CC$ using the marking
$\kappa$. Write $D=\{x_1,\ldots,x_r\}$ with $x_1<x_2<\ldots<x_r$ and
choose a base point $x_0\in U$ lying in the upper half plane. There is
a standard presentation
\begin{equation} \label{pi1pres}
   \pi_1(U,x_0) \;=\; \gen{\;\alpha_1,\ldots,\alpha_r \,\mid\,
                                  \prod_i\alpha_i = 1\;}
\end{equation}
of the fundamental group of $U$, depending only on $\kappa$. The
generators $\alpha_i$ is generated by a simple closed loop which
intersects the real line exactly twice, first the interval
$(x_{i-1},x_i)$, then the interval $(x_i,x_{i+1})$. 

Let $\V$ be a local system of $R$-modules on $U$, corresponding to a
representation $\rho:\pi_1(U,x_0)\to\GL(V)$. For $i=1,\dots,r$, set
$g_i:=\rho(\alpha_i)\in\GL(V)$.  Then we have
\[
      \prod_{i=1}^r \; g_i \;=\; 1.
\]
Conversely, given a free $R$-module $V$ of finite rank and a tuple
$\g=(g_1,\ldots,g_r)$ of elements of $\GL(V)$ satisfying the above relation,
we obtain a local system $\V$ which induces the tuple $\g$, as above.

\subsection{} \label{locals2}

We continue with the notation introduced in the previous subsection.
Let $j:U\inj X$ denote the inclusion. The {\em parabolic cohomology}
of $\V$ is defined as the sheaf cohomology of $j_*\V$, and is written as
\[
      H^n\para(U,\V) \;:=\; H^n(X,j_*\V).
\]
We have natural morphisms $H_c^n(U,\V)\to H\para^n(U,\V)$  and
$H\para^n(U,\V)\to H^n(U,\V)$ ($H_c$ denotes cohomology with compact
support). 

\begin{prop}  \label{locals2prop}
\begin{enumerate}
\item
  The group $H^n(U,\V)$ is canonically isomorphic to the group cohomology
  $H^n(\pi_1(U,x_0),V)$. In particular, we have
  \[
      H^0(U,\V) \;\cong\; V^{\gen{g_1,\ldots,g_r}},
  \]
  and $H^n(U,\V)=0$ for $n>1$. 
\item
  The map $H_c^1(U,\V)\to H\para^1(U,\V)$ is surjective and the map
  $H\para^1(U,\V)\to H^1(U,\V)$ is injective. In other words, $H\para^1(U,\V)$
  is the image of the cohomology with compact support in $H^1(U,\V)$. 
\end{enumerate}
\end{prop}

\proof Part (i) follows from the Hochschild--Serre spectral sequence and the
fact that the universal cover of $U$ is contractible. For (ii), see e.g.\ 
\cite{Looijenga92}, Lemma 5.3.  \Endproof

Let $\delta:\pi_1(U)\to V$ be a cocycle, i.e.\ we have
$\delta(\alpha\beta)=\delta(\alpha)\cdot\rho(\beta)+\delta(\beta)$. Set
$v_i:=\delta(\alpha_i)$. It is clear that the tuple $(v_i)$ is subject to the
relation
\begin{equation} \label{vcocyclerel}
   v_1\cdot g_2\cdots g_r\,+\,v_2\cdot g_3\cdots g_r\,+\ldots+\,v_r \;=\; 0.
\end{equation}
Conversely, any tuple $(v_i)$ satisfying \eqref{vcocyclerel} gives rise to a
unique cocycle $\delta$. This cocycle is a coboundary if and only if there
exists $v\in V$ such that $v_i=v\cdot(g_i-1)$ for all $i$. By Proposition
\ref{locals2prop} there is a natural inclusion
\[
      H^1\para(U,\V) \;\inj\; H^1(\pi_1(U),V).
\]
We say that $\delta$ is a {\em parabolic} cocycle if the class of $\delta$ in
$H^1(\pi_1(U),V)$ lies in the image of $H^1\para(U,\V)$.

\begin{lem}  \label{locals2lem}
  The cocycle $\delta$ is parabolic if and only if $v_i$ lies in the image of
  $g_i-1$, for all $i$.
\end{lem}

\proof Let $U_i\subset X$ be pairwise disjoint disks with center $x_i$, and set
$U_i^*:=U_i-\{x_i\}$. We have a long exact sequence
\begin{equation} \label{locals2lemeq2} 
   \cdots\to\; H^n_{x_i}(U_i,(j_!\V)|_{U_i}) \;\to\; H^n(U_i,(j_!\V)|_{U_i})
       \;\to\; H^n(U_i^*,\V|_{U_i^*}) \;\to\cdots.
\end{equation}
Given a class $c$ in $H^n(U_i,(j_!\V)|_{U_i})$ we can find a smaller
disk $U_i'\subset U_i$ with center $x_i$ such that the restriction of
$c$ to $U_i'$ vanishes (one way to see this is to use $\rm
\check{C}$ech cohomology). On the other hand, the cohomology groups
$H^n_{x_i}(U_i,(j_!\V)|_{U_i})$ and $H^n(U_i^*,\V|_{U_i^*})$ do not
change if we shrink the disk $U_i$. Therefore, by the exactness of
\eqref{locals2lemeq2} we have $H^n(U_i,(j_!\V)|_{U_i})=0$ and hence
\begin{equation} \label{locals2lemeq3}
   H^n_{x_i}(U_i,(j_!\V)|_{U_i}) \;\cong\; H^{n-1}(U_i^*,\V|_{U_i^*}) 
    \;\cong\;\;\begin{cases} 
                  \;{\rm Ker}(g_i-1),   & n=1,\\
                  \;{\rm Coker}(g_i-1), & n=2,\\
                  \;     0,             & \text{\rm otherwise}.
                \end{cases}
\end{equation}
For the second isomorphism we have used
$H^{n-1}(U_i^*,\V|_{U_i^*})\cong H^{n-1}(\pi_1(U_i),V)$. Consider
the long exact sequence 
\begin{equation} \label{locals2lemeq4}
      \cdots\to\; H^n_D(X,j_!\V) \;\to\; H^n_c(U,\V) \;\to\;
           H^n(U,\V) \;\to\cdots.
\end{equation}
By \eqref{locals2lemeq3}, the image in $H^2_D(X,j_!\V)$ of the class of a
$1$-cocycle $\delta:\pi_1(U)\to V$ vanishes if and only
$v_i:=\delta(\alpha_i)\in{\rm Im}\,(g_i-1)$. The lemma follows now from the
exactness of \eqref{locals2lemeq4} and from Proposition \ref{locals2prop}
(ii).  \Endproof

The preceding lemma shows that the association $\delta\mapsto(v_i)$ yields an
isomorphism
\[
   H^1\para(U,\V) \;\cong\; W_\g:=H_\g/E_\g,
\]
where
\[
   H_{\g} \;:=\; \{\,(v_1,\ldots,v_r) \;\mid\;
     v_i\in{\rm Im}(g_i-1),\; \text{\rm relation \eqref{vcocyclerel} holds}\;\}
\]
and
\[
   E_{\g} \;:=\; \{\, (\,v\cdot(g_1-1),\ldots,v\cdot(g_r-1)\,) \;\mid\;
                                    v\in V\;\}.
\]
The $R$-module $W_\g$ has already been defined in \cite{Voelklein01}, where it
is called the {\em braid companion} of $V$.

\begin{rem}  \label{locals2rem}
  Suppose that $R=K$ is a field and that the stabilizer $V^{\pi_1(U)}$ is
  trivial. Then the Ogg-Shafarevic formula implies the following dimension
  formula:
  \[
     \dim_K H^1\para(U,\V) \;=\; 
        (r-2)\dim_K V - \sum_{i=1}^r \dim_K{\rm Ker}(g_i-1).
  \] 
  This formula can also be verified directly using the isomorphism
  $H^1\para(U,\V)\cong W_\g$. 
\end{rem}

%-------------------------------------------------------------------

\section{Variation of a local system} \label{variation}
 
We study variations of local systems on the punctured sphere, with moving
singularities. The main result is the computation of the monodromy of the
parabolic cohomology of the variation. This computation is based on a
natural generalization of results of V\"olklein
\cite{Voelklein} \cite{Voelklein01}.

\subsection{}  \label{locals3}

Let $S$ be a connected complex manifold, and $r\geq 3$. An {\em
$r$-configuration} over $S$ consists of a smooth and proper morphism
$\pib:X\to S$ of complex manifolds together with a smooth relative
divisor $D\subset X$ such that the following holds. For all $s\in S$
the fiber $X_s:=\pib^{-1}(s)$ is a Riemann surface of genus $0$, and
the divisor $D\cap X_s$ consists of $r$ pairwise distinct points
$x_1,\ldots,x_r$.

Let us fix an $r$-configuration $(X,D)$ over $S$.  We set $U:=X-D$ and
denote by $j:U\inj X$ the natural inclusion.  Also, we write
$\pi:U\to S$ for the natural projection.  Choose a base point $s_0\in
S$ and set $X_0:=\pib^{-1}(s_0)$ and $D_0:=X_0\cap D$. Write
$D_0=\{x_1,\ldots,x_r\}$ and $U_0:=X_0-D_0=\pi^{-1}(s_0)$. Choose a
base point $x_0\in U_0$. The projection $\pi:U\to S$ is a topological
fibration and yields a short exact sequence
\begin{equation} \label{fibses}
  1 \;\To\; \pi_1(U_0,x_0) \;\To\; \pi_1(U,x_0) \;\To\; \pi_1(S,s_0) 
      \;\To\; 1.
\end{equation}

From now on, we shall drop the base points from our notation.
Let $\V_0$ be a local system of $R$-modules on $U_0$, corresponding to a
representation $\rho_0:\pi_1(U_0)\to\GL(V)$, as in \S \ref{locals1}. 

\begin{defn} \label{vardef}
  A {\em variation} of $\V_0$ over $S$ is a local system $\V$ of
  $R$-modules on $U$ whose restriction to $U_0$ is identified with
  $\V_0$. The {\em parabolic cohomology} of a variation $\V$ is the
  higher direct image sheaf
  \[
            \W \;:=\; R^1\pib_*(j_*\V).
  \]
\end{defn}

A variation $\V$ of $\V_0$ corresponds to a representation
$\rho:\pi_1(U)\to\GL(V)$ whose restriction to $\pi_1(U_0)$ is equal to
$\rho_0$. By definition, the parabolic cohomology
$\W$ of the variation $\V$ is a sheaf of $R$-modules on $S$. Locally
on $S$, 
the configuration $(X,D)$ is topologically trivial, i.e.\ there exists
a homeomorphism $X\cong X_0\times S$ which maps $D$ to $D_0\times
S$. It follows immediately that $\W$ is a local system with fibre
\[
     W \;:=\;  H^1\para(U_0,\V_0).
\]
In other words, $\W$ corresponds to a representation
$\eta:\pi_1(S)\to\GL(W)$. The following lemma provides a description
of $\eta$ in terms of cocycles.

\begin{lem}  \label{etalem}
  Let $\beta\in\pi_1(S)$ and $\delta:\pi_1(U_0)\to V$ be a parabolic cocycle.
  We write $[\delta]$ for the class of $\delta$ in $W$. Let
  $\betat\in\pi_1(U)$ be a lift of $\beta$. Then
  $[\,\delta\,]^{\eta(\beta)}=[\,\delta'\,]$, where $\delta':\pi_1(U_0)\to V$
  is the cocycle
  \[
      \alpha \;\longmapsto\; \delta(\betat\alpha\betat^{-1})\cdot
                                  \rho(\betat), \qquad
         \alpha\in\pi_1(U_0). 
  \]
\end{lem}

\proof We consider $\beta$ as a continuous map $\beta:I:=[0,1]\to S$.
Since $I$ is simply connected, there exists a continuous family of
homeomorphisms $\bar{\phi}_t:X_0\iso X_t:=\pib^{-1}(t)$, for $t\in
I$, such that $\bar{\phi}_t(D_0)=D_t:=X_t\cap D$ and such that
$\bar{\phi}_0$ is the identity. Let $\phi_t$ denote the restriction of
$\bar{\phi}_t$ to $U_0$. Note that $\phi_1:U_0\iso U_0$ is a
homeomorphism of $U_0$ with itself, whose homotopy class depends only
on $\beta\in\pi_1(S)$. We may further assume that
$\phi_1(x_0)=x_0$. Then $\betat:t\mapsto\phi_t(x_0)$ is a closed path
in $U$ with base point $x_0$. The class of $\betat$ in $\pi_1(U)$
(which we also denote by $\betat$) is a lift of $\beta\in\pi_1(S)$. It
is easy to check that for all $\alpha\in\pi_1(U_0)$ we have
\begin{equation} \label{etalemeq1}
          \phi_1(\alpha) \;=\; \betat^{-1}\,\alpha\,\betat.
\end{equation}
Since $\V$ is a local system on $U$, there exists a unique continuous
family $\psi_t:\V_0\iso\phi_t^*(\V|_{U_t})$ of isomorphisms of local
systems on $U_0$ such that $\psi_0$ is the identity on $\V_0$.
Evaluation of $\psi_t$ at the point $\betat(t)=\phi_t(x_0)$ yields a
continuous family of isomorphism $\psi_t(x_0):V\iso\V_{\betat(t)}$.
This family corresponds to a trivialization of $\betat^*\V$, and we
get
\begin{equation} \label{etalemeq2}
      \rho(\betat) \;=\; \psi_1(x_0).
\end{equation}
The pair $(\phi_t,\psi_t)$ induces a continuous family of
isomorphisms
\[
     \lambda_t: W \;\liso\; \W_t=H^1\para(U_t,\V|_{U_t}).
\]
Using \eqref{etalemeq1} and \eqref{etalemeq2}, one
finds that $\lambda_1([\,\delta\,])=[\,\delta'\,]$, where
\[
    \delta'(\alpha) \;=\; \psi_1(\delta(\phi_1^{-1}(\alpha)))
            \;=\; \delta(\betat\,\alpha\,\betat^{-1})\cdot\rho(\betat).
\]
By definition of the representation $\eta$, we have
$[\delta]^{\eta(\beta)}=\lambda_1([\delta])=[\delta']$. This completes
the proof of the lemma.  \Endproof

\begin{rem} \label{existencerem}

With the notation introduced above: let $\V_0$ be local system of
$R$-modules, corresponding to a representation
$\rho_0:\pi_1(U_0)\to\GL(V)$.  
\begin{enumerate}
\item
  A necessary condition for the existence of a variation of $\V_0$
  over $S$ is the following. For every element $\betat\in\pi_1(U,x_0)$
  there exists an element $g\in\GL(V)$ such that
  \[
      \rho_0(\betat\alpha\betat^{-1}) \;=\; g\rho_0(\alpha)g^{-1}
  \]
  holds for all $\alpha\in\pi_1(U_0)$. 
\item 
  Suppose that $S$ is a smooth affine curve, and that (i) holds. Then
  there exists a variation $\V$ of $\V_0$ over $S$. This follows
  easily from the fact that $\pi_1(S)$ is a free group. 
\item
  Suppose, moreover, that $R$ is an integral domain and that $\V_0$ is
  irreducible. If $\V'$ is another variation of $\V_0$ over $S$, then
  there exists a local system $\LL$ of rank one on $S$ such that
  $\V'\cong\V\otimes_R\pi^*\LL$. Let $\W$ (resp.\ $\W'$) denote the
  parabolic cohomology of $\V$ (resp.\ of $\V'$).  By the projection
  formula we have
  \[
      \W' \;\cong\; R^1\pib_*(j_*\V\otimes\pib^*\LL) \;\cong\;
           \W\otimes \LL.
  \]
  Therefore, the projective representation associated to $\W$,
  \[
       \lambda:\,\pi_1(S) \;\To\; \PGL(W),
  \]
  is uniquely determined by $\V_0$. 
\end{enumerate}
\end{rem}

\subsection{The Artin braid group and the cocycles $\Phi(\g,\beta)$} 
\label{locals5}
  
Let $D_0\subset\CC$ be a set of $r$ distinct complex numbers and set
$U_0:=\PP^1_\CC-D_0$. We choose a marking $\kappa$ of $(\PP^1_\CC,D_0)$ which
maps $\infty$ into the upper half plane, see \S \ref{locals1}. The choice of
$\kappa$ induces a presentation of $\pi_1(U_0,\infty)$, with generators
$\alpha_1,\ldots,\alpha_r$ and relation $\prod_i\alpha_i=1$.

Define
\[
       \OO_r \;:=\; \{\; D\subset \CC \;\mid\;
            |D|=r\;\}.
\]
The fundamental group $A_r \;:=\; \pi_1(\OO_r,D_0)$ is called the {\em
Artin braid group} on $r$ strands. The group $A_r$ has $r-1$
standard generators $\beta_1,\ldots,\beta_{r-1}$ with relations
\[
   \beta_i\beta_{i+1}\beta_i \;=\; \beta_{i+1}\beta_i\beta_{i+1},
     \qquad \beta_i\beta_j \;=\; \beta_j\beta_i,
\]
for $1\leq i<r$ and $i<j-1<r-1$. The element $\beta_i$ is represented by the
path $t\mapsto\{x_1,\ldots,\delta_i^-(t),\delta_i^+(t),\ldots,x_r\}$, where
$\delta_i^+$ (resp.\ $\delta_i^-$) is a path from $x_{i+1}$ to $x_i$ through
the inverse image under $\kappa$ of the upper half plane (resp.\ from $x_i$ to
$x_{i+1}$ through the inverse image of the lower half plane).

Define
\[
     \OO_{r,1} \;:=\; \{\;(D,x)\;\mid\; D\in\OO_r,\; x\in \PP^1_\CC-D \;\}.
\]
The natural projection $\OO_{r,1}\to\OO_r$ is a topological fibration with
fiber $U_0$, and admits a section $D\mapsto(D,\infty)$. It yields a
split exact sequence of fundamental groups
\begin{equation} \label{splitseq1}
   1 \;\to\; \pi_1(U_0,\infty) \;\To\; \pi_1(\OO_{r,1},(D_0,\infty)) 
     \;\To\; A_r \;\to\; 1.
\end{equation} 
We may identify $A_r$ with its image in $\pi_1(\OO_{r,1})$ under the
splitting induced from the section $D\mapsto(D,\infty)$. Then
$A_r$ acts, by conjugation, on $\pi_1(U_0,\infty)$. We have the
following well known formulas for this action:
\begin{equation} \label{locals5eq2}
  \beta_i^{-1}\,\alpha_j\beta_i \;=\quad
  \begin{cases}
    \quad \alpha_i\alpha_{i+1}\alpha_i^{-1}, & 
                              \quad\text{\rm for $j=i$,} \\
    \quad \alpha_i, & \quad \text{\rm for $j=i+1$,} \\
    \quad \alpha_j, & \quad \text{\rm otherwise.} 
  \end{cases}
\end{equation}

Let $R$ be a commutative ring and $V$ a free $R$-module of finite
rank. Define
\[
    \E_r \;:=\; \{\;\g=(g_1,\ldots,g_r) \mid 
                  g_i\in\GL(V),\;\; \prod_ig_i=1 \;\}.
\]
An element $\g\in\E_r$ corresponds to a representation
$\rho_0:\pi_1(U_0)\to\GL(V)$ (set $\rho_0(\alpha_i):=g_i$) and
hence to a local system $\V_0$ on $U_0$.  Given $\beta\in A_r$, we
set
\begin{equation} \label{locals5eq3}
    \rho_0^{\beta}(\alpha) \;:=\; \rho_0(\beta\alpha\beta^{-1})
\end{equation}
and call the local system $\V_0^{\beta}$ corresponding to the
representation $\rho_0^{\beta}$ the {\em twist} of $\V_0$ by
$\beta$. We denote by $\g^{\beta}$ the element of $\E_r$ corresponding
to $\rho_0^{\beta}$. This defines an action of $A_r$ on $\E_r$, from
the right. From \eqref{locals5eq2} we get the following formula for
the effect of the standard generators $\beta_i$ on $\E_r$:
\begin{equation}  \label{locals5eq4}
   \g^{\beta_i} \;=\; (g_1,\ldots,g_{i+1},g_{i+1}^{-1}g_ig_{i+1},\ldots,g_r).
\end{equation}

Given $\g\in\E_r$, we have defined in \S \ref{locals2} the $R$-module
\[
    H_\g \;=\; \{\;(v_1,\ldots,v_r) \mid v_i\in{\rm Im}(g_i-1),\;
            \text{\rm relation $\eqref{vcocyclerel}$ holds}\;\}.
\]
An element $(v_1,\ldots,v_r)\in H_{\g}$ corresponds to a parabolic
cocycle $\delta:\pi_1(U_0)\to V$, determined by
$\delta(\alpha_i)=v_i$, for $i=1,\ldots,r$. Here the
$\pi_1(U_0)$-module structure on $V$ is induced by $\g$. We say that
$\delta$ is a parabolic cocycle with respect to $\g$. Given $\beta\in
A_r$ and $\alpha\in\pi_1(U_0)$, set
\[
     \delta^{\beta}(\alpha) \;:=\; \delta(\beta\,\alpha\,\beta^{-1}).
\]
One easily checks that $\delta^{\beta}:\pi_1(U_0)\to V$ is a
parabolic cocycle with respect to $\g^{\beta}$. Moreover, the
association $\delta\mapsto\delta^{\beta}$ defines an $R$-linear map
\[
      \Phi(\g,\beta):\,H_{\g} \;\To\; H_{\g^{\beta}}.
\]
In order to maintain compatibility with our convention of `acting from
the right', we write $(v_i)^{\Phi(\g,\beta)}$ for the image of
$(v_i)\in H_{\g}$ under $\Phi(\g,\beta)$. Using \eqref{locals5eq2}
and the fact that $\delta$ is a cocycle with respect to $\g$, we get
\begin{multline} \label{locals5eq6}
 \qquad (v_1,\ldots,v_r)^{\Phi(\g,\beta_i)} \\ 
    \;=\;  (v_1,\ldots,v_{i+1},\,
    \underbrace{v_{i+1}(1-g_{i+1}^{-1}g_ig_{i+1})+v_ig_{i+1}}_{
        \text{\rm $(i+1)$th entry}},\,\ldots,v_r).
\end{multline} 
Moreover, we have the `cocycle rule'
\begin{equation} \label{locals5eq7}
    \Phi(\g,\beta)\cdot\Phi(\g^{\beta},\beta') \;=\; 
       \Phi(\g,\beta\beta').
\end{equation}
(The product on the left hand side of \eqref{locals5eq7} is defined as
the function from $H_{\g}$ to $H_{\g^{\beta\beta'}}$ obtained by
first applying $\Phi(\g,\beta)$ and then
$\Phi(\g^{\beta},\beta')$.) 

The submodule 
\[
      E_{\g} \;:=\; \{\;(v\cdot(g_1-1),\ldots,v\cdot(g_r-1) \;\mid\;
                           v\in V \;\}
\]
of $H_{\g}$ corresponds to cocycles $\delta$ which are
coboundaries. It is easy to see that $\Phi(\g,\beta)$ maps $E_{\g}$
into $E_{\g^{\beta}}$ and therefore induces an isomorphism
\[
       \bar{\Phi}(\g,\beta):\,W_{\g}:=H_{\g}/E_{\g} \;\liso\;
              W_{\g^{\beta}}.
\]
One can compute $\bar{\Phi}(\g,\beta)$ explicitly for all $\beta\in
A_r$ using \eqref{locals5eq6} and \eqref{locals5eq7}, provided that
$\beta$ is given as a word in the standard generators $\beta_i$.
Moreover, this computation can easily be implemented on a computer.

Given $\g\in\E_r$ and $h\in\GL(V)$, we set 
\[
         \g^h \;=\; (h^{-1}g_1h,\ldots,h^{-1}g_r h),
\]
and we define an isomorphism
\[
   \Psi(\g,h):\; \left\{\;
     \begin{array}{ccc}
        H_{\g\cdot h}  &  \;\liso\; & H_{\g} \\
        (v_1,\ldots,v_r) & \;\longmapsto\; & 
             (v_1\cdot h,\ldots,v_r\cdot h).
     \end{array}\right.
\]
It is clear that $\Psi(\g,h)$ maps $E_{\g\cdot h}$ to $E_\g$ and
therefore induces an isomorphism $\bar{\Psi}(\g,h):W_{\g\cdot h}\iso
W_\g$.

\subsection{Explicit computation of the monodromy}  \label{locals6}

Let us go back to the situation of \S \ref{locals3}: we are given an
$r$-configuration $(X,D)$ over a connected complex manifold $S$. We
have also chosen a base point $s_0\in S$. As usual, we set $U:=X-D$,
and denote by $U_0$ the fiber of $U\to S$ over $s_0$.

\begin{defn} \label{framedef}
  An {\em affine frame} for the configuration $(X,D)$ is an isomorphism of
  complex manifolds $\lambda:X\cong\PP^1_S$, compatible with the
  projection to $S$, such that $\lambda(D)$ either contains or is
  disjoint from $\{\infty\}\times S$. 
\end{defn}

In this subsection, we shall assume that there exists an affine frame for
$(X,D)$, and we use it to identify $X$ with $\PP^1_S$. We remark that there
exist configurations $(X,D)$ which do not admit an affine frame (e.g.\ because
$X\not\cong\PP^1_S$). It seems, however, that such examples have no
practical relevance for the problems this paper is about. 

By the nature of Definition \ref{framedef}, there are two cases to
consider.  Suppose first that $D$ is disjoint from $\{\infty\}\times
S$. Then $D$ gives rise to a map $p:S\to\OO_r$ which sends $s\in S$ to
the fiber of $D\subset\AA^1_S\to S$ over $s$.  Set
$D_0:=p(s_0)\subset\CC$. Choose a marking $\kappa$ of
$(\PP^1_\CC,D_0)$. We will use $\kappa$ to identify the fundamental
group $\pi_1(\OO_r,D_0)$ with the Artin braid group $A_r$, as in the
previous subsection.  Let $\varphi:\pi_1(S,s_0)\to A_r$ denote the
group homomorphism induced by $p$.  The exact sequence
\begin{equation} \label{splitseq2}
   1 \;\to\; \pi_1(U_0,\infty) \;\To\; \pi_1(U,(\infty,s_0)) 
     \;\To\; \pi_1(S,s_0) \;\to\; 1
\end{equation} 
of the fibration $U\to S$ can be identified with the
pullback of the sequence \eqref{splitseq1} along $\varphi$. Using the
splitting of \eqref{splitseq2} coming from the $\infty$-section, we
will consider $\pi_1(S)$ as a subgroup of $\pi_1(U)$. By construction,
the action of $\pi_1(S)$ on $\pi_1(U_0,\infty)$ factors through the
map $\varphi$ and is given by the formulas \eqref{locals5eq2}. 

Now suppose that $D$ contains the section $\{\infty\}\times S$. We
denote by $\pi_1(U_0,\infty)$ the fundamental group of $U_0$ with
$\infty$ as `tangential base point'. More precisely, consider subsets
of $U_0\subset\CC$ of the form $\Omega_t=\{z\in\CC\mid
|z|>t,\,z\not\in(-\infty,0)\,\}$, for $t\gg 0$. The fundamental group
$\pi_1(U_0,\Omega_t)$ is independent of $t$, up to {\em canonical}
isomorphism, so we may define
$\pi_1(U_0,\infty):=\lim\pi_1(U_0,\Omega_t)$. With this convention,
the sequence \eqref{splitseq2} is still well defined and admits a
canonical section. In fact, the fibration $U\to S$ admits a section
$\xi:S\to U$, unique up to homotopy, such that for all $s\in S$ we
have $\xi(s)\in\Omega_t\subset U_s$, for some $t>0$. As in the first
case, we will identify $\pi_1(S)$ with the image of this section.

The {\em Hurwitz braid group} $B_r$ is defined as the fundamental
group of the set
\[
   \U_r \;=\; \{\,D\subset\PP^1_\CC\, \mid\, |D|=r \,\},
\]
with base point $D_0$.  The natural map $\OO_r\to\U_r$ identifies
$B_r$ with the quotient of $A_r$ by the relation
\[
     \beta_1\beta_2\cdots\beta_{r-1}^2\cdots\beta_2\beta_1 \;=\; 1.
\]
The configuration $(X,D)$ induces a map $p:S\to\U_r$ and a
homomorphism $\varphi:\pi_1(S)\to B_r$. If $\{\infty\}\times S\subset
D$ then the image of $\varphi$ is contained in the subgroup of $B_r$
generated by the first $r-2$ standard braids
$\beta_1,\ldots,\beta_{r-2}$, which is isomorphic to $A_{r-1}$.
Moreover, just as in the first case, the action of $\pi_1(S,s_0)$ on
$\pi_1(U_0,\infty)$ by conjugation factors through the map $\varphi$
and is given by the formulas \eqref{locals5eq2}. From now on, we will
treat both cases of Definition \ref{framedef} simultaneously.

Let $\V_0$ be a local system of free $R$-modules on
$U_0=\PP^1(\CC)-D_0$, corresponding to a representation
$\rho_0:\pi_1(U_0)\to\GL(V)$. A variation of $\V_0$ over $S$
corresponds, by definition, to a representation
$\rho:\pi_1(U)\to\GL(V)$ whose restriction to $\pi_1(U_0)$ equals
$\rho_0$. Obviously, $\rho$ is uniquely determined by its restriction
to $\pi_1(S)$, which we denote by $\chi:\pi_1(S)\to\GL(V)$. Then
\begin{equation} \label{chieq}
       \rho_0(\gamma\,\alpha\,\gamma^{-1}) \;=\; 
           \chi(\gamma)\,\rho_0(\alpha)\,\chi(\gamma)^{-1}
\end{equation}
holds for all $\alpha\in\pi_1(U_0)$ and $\gamma\in\pi_1(S)$.  With
$\g\in\E_r$ corresponding to $\rho_0$ (via the choice of the marking
$\kappa$), this is equivalent to
\begin{equation} \label{locals6eq2}
         \g^{\varphi(\gamma)} \;=\; \g^{\chi(\gamma)^{-1}}.
\end{equation}
Let $\W$ be the parabolic cohomology of $\V$ and
$\eta:\pi_1(S)\to\GL(W_{\g})$ the corresponding representation (here
we identify the fiber of $\W$ at $s_0$ with the $R$-module
$W_{\g}=H_{\g}/E_{\g}$, see the previous subsection).

\begin{thm} \label{etathm}
  For all $\gamma\in\pi_1(S)$ we have
  \[
      \eta(\gamma)  \;=\;
           \bar{\Phi}(\g,\varphi(\gamma))\cdot\bar{\Psi}(\g,\chi(\gamma)),
  \]
  where $\bar{\Phi}(\g,\beta):W_{\g}\iso W_{\g^{\beta}}$ and
  $\bar{\Psi}(\g,h):W_{\g^h}\to W_{\g}$ are the isomorphisms defined
  in \S \ref{locals5}.
\end{thm}

\proof
Straightforward, using Lemma \ref{etalem}, the definition of
$\bar{\Phi}(\g,\beta)$ and $\bar{\Psi}(\g,h)$, and \eqref{locals6eq2}.
\Endproof

%----------------------------------------------------------------

\section{\'Etale local systems} \label{etale}

We transfer the situation considered in the first two sections into the
\'etale world, and we state a comparison theorem. We also prove a theorem
which is useful to bound the field of linear moduli of an (\'etale) local
system which is obtained as the parabolic cohomology of a variation.

\subsection{Recall} \label{etale0}

In this section, we fix a prime number $l$ and a finite extension
$K/\QQ_l$. We denote by $R$ one of the following rings:
(a) $R:=K$, (b) $R:=\OO_K$, the ring of integers of $K$, or
(c) $R:=\OO_K/\ell^m$, where $\ell$ is the prime ideal of $\OO_K$.

Let $k$ be a field of characteristic $0$ and $X$ a smooth,
geometrically irreducible scheme over $k$. Also, let
$x:\Spec\bar{k}\to S$ be a geometric point.  We denote by
$\pi_1(X)=\pi_1(X,x)$ the algebraic fundamental group of $X$ with base
point $x$.

An {\em \'etale local system} $\V$ of $R$-modules on $X$ is, by
definition, a {\em locally constant and constructible sheaf of
$R$-modules} \cite{MilneEC}, whose stalks are free $R$-modules of
finite rank. In case $R=K$, this is also called a {\em lisse $\ell$-adic
sheaf} \cite{SGA4}. It is a standard fact that $\V$ corresponds to a
continuous representation
\[
      \rho:\,\pi_1(X,x) \;\To\; \GL(V),
\]
where $V:=\V_x$ is the stalk of $\V$ at $x$.  

Now suppose that $k\subset\CC$ is a subfield of the complex numbers.  The set
of $\CC$-rational points of $X$ has a canonical structure of a complex
manifold, which we denote by $X^{\rm an}$. Moreover, there is a functor
$\F\mapsto\F\an$ from sheaves (of abelian groups) on $X_{\et}$ to sheaves on
$X\an$, called {\em analytification} (see e.g.\ \cite{FreitagKiehl}, \S I.11).
If $\V$ is an \'etale local system on $X$ corresponding to a representation
$\rho:\pi_1(X,x)\to\GL(V)$, then the analytification $\V\an$ of $\V$ is the
local system corresponding to the composition of $\rho$ with the natural
homomorphism $\pi_1\tp(X^{\rm an},x)\to\pi_1(X,x)$.

\subsection{Parabolic cohomology of an \'etale local system}  \label{etale1}

Let $k$ be a field of characteristic $0$ and let $S$ be a
smooth, affine and geometrically connected variety over $k$. Let
$(X,D)$ be an {\em $r$-configuration} over $S$. By this we mean that
$\pib:X\to S$ is a proper smooth curve of genus $0$ and $D\subset X$
is a smooth relative divisor of relative degree $r$ (compare with \S
\ref{locals3}).  We denote by $j:U:=X-D\inj X$ the inclusion and by
$\pi:U\to S$ the natural projections. We 
fix a $k$-rational point $s_0$ on $S$ as a base point. We write
$U_0:=\pi^{-1}(s_0)$ for the fiber of $\pi$ over $s_0$ and we
choose a geometric point $x_0:\Spec\kb\to U_0$ as base point.

\begin{defn} \label{varconvdef}
  Let $\V_0$ be an \'etale local system of $R$-modules on
  $U_0$. A {\em variation} of $\V_0$ over $S$ is an \'etale
  local system $\V$ on $U$ whose restriction to $U_0$ is equal
  to $\V_0$. The {\em parabolic cohomology} of the variation $\V$ is
  the sheaf of $R$-modules on $S_{\acute{e}t}$
  \[ 
        \W \;:=\; R^1\pib_*(j_*\V).
  \]
  See \cite{MilneEC}. 
\end{defn}

\begin{thm} \label{etalethm}
  Suppose that $k\subset\CC$. 
  \begin{enumerate}
  \item
    $\W$ is an \'etale local system of $R$-modules. 
  \item 
    There is a natural isomorphism of local systems of $R$-modules on $X\an$
    \[
         \W\an \;\liso\; R^1\pib_*(j_*\V\an),
    \]
    functorial in $\V$.
  \end{enumerate}
\end{thm}

\proof Using standard arguments (see e.g.\ \cite{FreitagKiehl}, \S 12), one
reduces the claim to the case $R=\OO_K/\ell^m$. Let $\F$ be a constructible
sheaf of $R$-modules on $X$. By the comparison theorem between \'etale and
singular cohomology there is a natural isomorphism of sheaves on $S\an$
\[
    (R^1\pib_*\F)\an \;\liso\; R^1\pib_*(\F\an).
\]
(See e.g.\ \cite{FreitagKiehl}, Theorem 11.6, for the case where $k=\CC$. The
general case follows immediately, using the Proper Base Change Theorem,
\cite{FreitagKiehl}, Theorem 6.1.) It is easy to see (e.g.\ using
\cite{FreitagKiehl}, Proposition 11.4) that $(j_*\V)\an=j_*(\V\an)$.
Therefore, Part (ii) of the theorem follows from the comparison theorem.
By \cite{FreitagKiehl}, Theorem 8.10, the sheaf $\W=R^1\pib_*(j_*\V)$ is
constructible. But we have just proved that $\W\an$ is locally constant, which
shows that $\W$ is locally constant as well.  This finishes the proof of the
theorem.  \Endproof

\subsection{The field of linear moduli}  \label{etale3}

As in the previous subsection, $S$ denotes a smooth and geometrically
connected $k$-variety and $(X,D)$ an $r$-configuration over $S$. 
We assume that $k\subset\CC$ and denote by $\kb$ the algebraic closure of $k$
inside $\CC$. 

Unlike in the previous subsection, $\V_0$ is now an \'etale local system of
$R$-modules on the {\em geometric} fibre $U_{0,\kb}:=U_0\otimes\kb$, and $\V$
denotes a variation of $\V_0$ over $S_{\kb}:=S\otimes_k\kb$. Let $\W$ be the
parabolic cohomology of $\V$. By construction, $\W$ is an \'etale local system
of $R$-modules on $S_{\kb}$.  For $\sigma\in\Gal(\kb/k)$, we denote by
$\tilde{\sigma}_{U_0}:U_{0,\kb}\iso U_{0,\kb}$ the semi-linear automorphism
corresponding to $\sigma$ and the $k$-model $U_0$. The {\em twist of $\V_0$ by
  $\sigma$} (with respect to the $k$-model $U_0$) is the \'etale local system
$\V_0^\sigma:=\tilde{\sigma}_{U_0}^*\V_0$.

\begin{defn} \label{folmdef}
  We say that $k$ is a {\em field of linear moduli} for $\V_0$ if the \'etale
  local system $\V_0^\sigma$ is isomorphic to $\V_0$, for all
  $\sigma\in\Gal(\kb/k)$. We say that $k$ is a {\em field of projective
    moduli} for $\V_0$ if for all $\sigma\in\Gal(\kb/k)$ there exists an
  \'etale local system $\LL_\sigma$ of rank one such that
  $\V_0^\sigma\cong\V_0\otimes_R\LL_\sigma$. Similarly, one defines the notion
  of `field of linear/projective moduli' for the local systems $\V$ and $\W$. 
\end{defn}

\begin{thm} \label{folmthm}
\begin{enumerate}
\item
  If $k$ is a field of linear moduli for $\V$ then it is also a field of linear
  moduli for $\W$. 
\item Suppose that $R=\OO_K$ or $R=K$, and that $\V_0$ is irreducible. Then if
  $k$ is a field of projective moduli for $\V_0$, it is also a field of
  projective moduli for $\W$. Moreover, the projective representation
  $\lambda\geo:\pi_1(S_{\kb})\to\PGL(W)$ associated to $\W$ extends to a
  projective representation $\lambda:\pi_1(S)\to\PGL(W)$.
\end{enumerate}
\end{thm}
 
\proof Let $\sigma\in\Gamma_{k_0}$ with $\V^\sigma\cong\V$.  Using the Proper
Base Change Theorem we get
\[
   \W^\sigma \;=\; \hat{\sigma}_S^* R^1\pib_{\kb,*}(j_{\kb,*}\V) \;=\;
     R^1\pib_{\kb,*}(j_{\kb,*}\V^\sigma) \;\cong\; \W.
\]
This proves (i).  The proof of (ii) is a combination of the preceeding
argument and Remark \ref{existencerem} (iii).  \Endproof

\begin{rem} \label{folmrem}
  Theorem \ref{folmthm} can be used to give a new proof of the main
  result of \cite{VoelkleinCrelle} (which essentially states that the
  braid companion functor preserves the field of linear moduli). In
  the rest of the paper, we will not use Theorem \ref{folmthm} and we
  will not need the concept `field of linear moduli'. The point is
  that in our main example in \S \ref{example}, the variation $\V$ is
  already known to be defined over $\QQ$, which means that the
  resulting local system $\W$ is defined over $\QQ$ as well, by
  construction. In fact, this seems to be the case for all known
  applications of these and similar methods to the Regular Inverse
  Galois Problem. Nevertheless, the authors think that Theorem
  \ref{folmthm} may be useful for future applications.
\end{rem}

%---------------------------------------------------------------------

\section{Local systems on Hurwitz curves}  \label{hurwitz}

A finite Galois cover $f:Y\to\PP^1$ together with a representation
$G\inj\GL_n(K)$ of its Galois group corresponds to a local system on
$\PP^1$ with finite monodromy. Therefore, a representation
$G\inj\GL_n(K)$ gives rise to a variation of local systems on a
certain Hurwitz space $H$. Since Hurwitz spaces are algebraic
varieties, the parabolic cohomology of this variation corresponds to a
Galois representations of the function field of $H$. In case $H$ is a
rational variety, this has potential applications to the Regular
Inverse Galois Problem.

\subsection{} \label{hurwitz1}

In this section we fix a finite group $G$. For each integer $r\geq 3$, we set
\[
    \E_r(G) \;:=\; \{\,\g=(g_1,\ldots,g_r)\;\mid\;
       G=\gen{g_i},\;\prod_ig_i=1\;\}.
\]
An element $\g$ of this set is called a {\em generating system} of
length $r$ for the group $G$. The group $G$ acts on the set $\E_r(G)$
by simultaneous conjugation. We write $\Ni_r(G)$ for the sets of
orbits of this action. Elements of $\Ni_r(G)$ are called {\em Nielsen
classes} and written as $[\g]$, with $\g\in\E_r(G)$.

The Artin braid group $A_r$ acts on the set $\E_r(G)$ from the right,
in a standard way, see e.g.\ \cite{FriedVoe91} and \S \ref{locals5}.
This action extends to an action of its quotient $B_r$, the Hurwitz
braid group. By abuse of notation, we denote the image of the standard
generator $\beta_i\in A_r$ in $B_r$ by the same name.

Suppose, for the moment, that $r=4$. The elements $\beta_1\beta_3^{-1}$ and
$(\beta_1\beta_2\beta_3)^2$ generate a normal subgroup $Q\lhd B_4$,
isomorphic to the Klein $4$-group. The quotient $\Bb_4:=B_4/Q$ is called the
{\em mapping class group}. The set of $Q$-orbits of $\Ni_r(G)$ is denoted by
$\Ni_r\red(G)$. Elements of this set are called {\em reduced Nielsen classes},
and are written as $[\g]\red$. The action of $B_4$ on $\Ni_4(G)$
descents to an action of the mapping class group $\Bb_4$ on $\Ni_4\red(G)$.

Let $\bC=(C_1,\ldots,C_r)$ be an ordered $r$-tuple of conjugacy
classes of the group $G$. We say that $\g\in\E_r(G)$ has {\em type
$\bC$} if there exist an integer $n$, prime to the order of $G$, and a
permutation $\sigma\in S_r$ such that $g_i^n\in C_{\sigma(i)}$ for
$i=1,\ldots,r$. The subset of $\E_r(G)$ of all elements of type $\bC$
is denoted by $\E(\bC)$. We also obtain subsets
$\Ni(\bC)\subset\Ni_r(G)$ and $\Ni\red(\bC)\subset\Ni_4\red(G)$.

\subsection{} \label{hurwitz2}

Let $T$ be a scheme over $\Spec(\QQ)$ and $X$ a smooth projective
curve over $T$ of genus $0$. A {\em $G$-cover} of $X$ is a finite
morphism $f:Y\to X$, together with an isomorphism $G\cong{\rm
Aut}(Y/X)$, such that the following holds. First, $f$ is tamely
ramified along a smooth relative divisor $D\subset X$ with constant
degree $r:={\rm deg}(D/S)$, and \'etale over $U:=X-D$. Second, for
each geometric point $t:\Spec k\to T$, the pullback $f_t:Y_t\to X_t$
is a $G$-Galois cover (in particular, $Y_t$ is connected). We say that
two $G$-covers $f_1:Y_1\to X_1$ and $f_2:Y_2\to X_2$ defined over $T$
are {\em isomorphic} if there exist isomorphisms of $T$-schemes
$\psi:Y_1\iso Y_2$ and $\phi:X_1\iso X_2$ such that $\phi\circ
f_1=f_2\circ\psi$.

By the result of \cite{FriedVoe91}, \cite{DebesFried} and
\cite{diss}, there exists a certain $\QQ$-scheme, denoted by
\[
     H\red_r(G),
\]
which is a coarse moduli space for $G$-Galois covers of curves of genus $0$
with $r$ branch points. In particular, to each $G$-Galois cover $f:Y\to X$
over a $\QQ$-scheme $T$, we can associate a map $\varphi_f:T\to H$, called the
{\em classifying map} for $f$. The association $f\mapsto\varphi_f$ is
functorial in $T$. For $T=\Spec k$, where $k$ is an algebraically closed
field of characteristic $0$, it induces a bijection between isomorphism
classes of $G$-covers of curves of genus $0$ with $4$ branch points, defined
over $k$, and $k$-rational points on $H$. The scheme $H\red_r(G)$ is called
the {\em reduced inner Hurwitz space}. It is a smooth and affine variety over
$\QQ$. 

For the rest of this section, we will assume that $r=4$. In this case,
the variety $H\red_4(G)$ is a smooth affine curve,
equipped with a finite separable cover
\[
        j:H_4\red(G) \;\To\; \AA^1_{\QQ},
\]
which is at most tamely ramified at $0,1728$ and \'etale over
$\AA^1_{\QQ}-\{0,1728\}$.  The map $j$ is characterized by the
following property. Let $t:\Spec k\to H_4\red(\bC)$ be a geometric
point, corresponding to a $G$-cover $f:Y\to X$ with branch locus
$D=\{x_1,\ldots,x_4\}$. Then $j(t)\in k$ is the $j$-invariant of the
configuration $(X,D)$.

Let $x_1,\ldots,x_4\in\CC$ be four distinct complex numbers. Set
$X_0:=\PP^1_\CC$ and $D_0:=\{x_1,\ldots,x_4\}\subset X_0$. Let $z_0\in\RR$
denote the $j$-invariant of the configuration $(\PP^1_\CC,D_0)$. We assume
that $z_0\not=0,1728$. After choosing a marking $\kappa$ of $(X_0,D_0)$, we
obtain a presentation of $\pi_1(U_0)$ with generators
$\alpha_1,\ldots,\alpha_4$ and relation $\prod_i\alpha_i=1$. This presentation
yields a bijection
\begin{equation} \label{Nieq}
   j^{-1}(z_0) \;\liso\; \Ni_4\red(G).
\end{equation}
(Recall that the left hand side of \eqref{Nieq} may be identified with
the set of isomorphism classes of $G$-covers of $X_0$ with branch
locus $D_0$.) The fundamental group of $\CC-\{0,1728\}$ acts on the
left hand side of \eqref{Nieq}, and the mapping class group $\Bb_4$
acts on the right hand side. There is a natural identification of
these two groups which makes the bijection \eqref{Nieq}
equivariant. In particular, we obtain a bijection between the set of
connected components of $H_4\red(G)_\CC$ and the $\Bb_4$-orbits of
$\Ni_4\red(G)$.

\subsection{} \label{hurwitz3}

Let us now fix the following objects:
\begin{itemize}
\item
  an $\Bb_4$-orbit $O\subset\Ni_4\red(G)$, and 
\item 
  a faithful and irreducible linear representation $G\inj\GL_n(K)$, with
  coefficients in a number field $K$.
\end{itemize}
The orbit $O$ corresponds to a connected component $H(O)$ of
$H_4\red(G)_\CC$. Let $k_O$ be its field of definition, i.e.\ the
smallest subfield $k\subset\CC$ such that the natural action of
$\Aut(\CC/k)$ on the set of connected components of $H_4\red(G)$
stabilizes $H(O)$. We may and will consider $H(O)$ as a geometrically
irreducible variety over $k_O$ (i.e.\ as a connected component
$H_4\red(G)\otimes k_O$) .

Let $f:Y\to X$ be a {\em versal family} over $H(O)$, i.e.\ a $G$-cover
defined over a $k_O$-scheme $S$ whose classifying map $S\to
H_4\red(G)$ is \'etale and factors through the natural map $H(O)\to
H_4\red(G)$. (It is known that such a versal family always
exists. Under some favorable conditions, one can take $S=H(O)$, see
\S \ref{example1}.) Let $k$ be the field of definition of $S$, i.e.\
the algebraic closure of $\QQ$ inside the function field of $S$. Since
the map $S\to H\red(G)$ is \'etale, $S$ is a smooth, affine and
geometrically connected $k$-curve. Moreover $k$ is a finite extension
of $k_O$.

Choose an integer $N$ such that $G\subset\GL_n(K)$ is contained in
$\GL_n(R)$, where $R=\OO_K[1/N]$.  Let $D\subset X$ be the branch locus
of $f:Y\to X$, and set $U:=X-D$. The $G$-cover $f:Y\to X$ gives rise to a
surjective group homomorphism $\pi_1(U)\to G$. We denote by
\[
    \rho:\,\pi_1(U)\;\To\;\GL_n(R)
\]
the composition of this homomorphism with the injection $G\inj\GL_n(R)$.

We shall write $\rho\geo:\pi_1(U_{\Qb})\to\GL_n(R)$ (resp.\ 
$\rho\tp:\pi_1\tp(U_\CC)\to\GL_n(R)$) for the restriction of $\rho$ to the
geometric fundamental group of $U$ (resp.\ to the topological fundamental
group of the analytic space associated to $U$). Let $\V\an$ denote the local
system of $R$-modules on $U_\CC$ corresponding to $\rho\tp$. Also, let 
\[
   \W\an \;:=\;  R^1\pib_*(j_*\V\an)
\]
be the parabolic cohomology of $\V\an$. Recall that $\W\an$ is a local
system of $R$-modules corresponding to a representation
$\eta\tp:\pi_1\tp(S_\CC)\to\GL(W)$.

On the other hand, for each prime ideal $\p$ of $K$ which is prime to
$N$ we let
\[
    \rho_\p:\pi_1(U) \;\To\; \GL_n(\OO_\p)
\]
denote the $p$-adic representation induced by $\rho$. It corresponds to an
\'etale local system $\V_\p$ of $\OO_\p$-modules on $U$. Again we can
form the (parabolic) higher direct image of $\V_\p$,
\[
    \W_\p \;:=\; R^1\pib_*(j_*\V_\p),
\]
which is an  \'etale local system of $\OO_\p$-modules, thus  corresponds to a
representation $\eta_\p:\pi_1(S)\to\GL(W_\p)$. 

\begin{prop} 
  There exists a canonical isomorphism $W_\p\cong
W\otimes_R\OO_\p$ such that the following diagram commutes:
\[\begin{CD}
  \pi_1\tp(S_\CC)  @>{\eta\tp}>>  \GL(W)     \\
  @VVV                                @VVV       \\
  \pi_1(S)   @>{\eta_\p}>> \GL(W_\p). \\
\end{CD}\]
Hence the image of $\eta_\p\geo$ is equal to the topological
closure of the image of $\eta\tp$.
\end{prop}

\proof This follows from Theorem \ref{etalethm} and the fact that
$\pi_1(S_{\Qb})$ is the profinite completion of $\pi_1\tp(S_\CC)$.
\Endproof

\subsection{}  \label{hurwitz4}

We can now use the results of \S \ref{locals6} to determine the image
of $\eta\tp$. We will use the notation introduced in the previous
subsection, with the following difference. Since we will be working
exclusively with complex analytic spaces, we will omit the index $(\ 
)_\CC$. For instance, we write $S$ instead of $S_\CC$, etc. 

Choose a point $s_0\in S$ with $j(s_0)=z_0$ and let $(X_0,D_0)$ denote
the fibre of the configuration $(X,D)$ over $s_0$. Let
$\rho_0\tp:\pi_1\tp(U_0,x_0)\to\GL_n(R)$ denote the restriction of
$\rho\tp$ to the subgroup $\pi_1\tp(U_0,x_0)\subset\pi_1\tp(U,x_0)$,
and set $g_i=\rho_0(\alpha_i)\in G$. By construction, the tuple
$\g:=(g_i)$ is a generating system for $G$, and the reduced Nielsen
class of $\g$ is an element of the $\Bb_4$-orbit $O$. Moreover,
$[\g]\red$ is stabilized by the image of the the group homomorphism 
\[
           \varphib:\pi_1\tp(S,s_0) \;\To\; \Bb_4,
\] 
which is induced by the configuration $(X,D)$ over $S$. 

For simplicity, we also assume that the configuration $(X,D)$ admits
an affine frame (Definition \ref{framedef}), which we use to identify
$X$ with $\PP^1_S$. This assumption will be satisfied in our main
example. Note also that (at least in the situation where $S$ is
one-dimensional), there always exists an affine frame over a dense
open subset of $S$.

The $\infty$-section defines a section of the natural projection
$\pi_1\tp(U)\to\pi_1\tp(S)$. We identify $\pi_1\tp(S)$ with the
image of this section. Let $\chi:\pi_1\tp(S)\to G$ denote the restriction
of $\rho$ to $\pi_1\tp(S)$. Essentially by definition, we have
\begin{equation} \label{ggeq}
        \g^{\varphi(\gamma)} \;=\; \g^{\chi(\gamma)^{-1}},
\end{equation}
for all $\gamma\in\pi_1(S)$, compare with
\eqref{locals6eq2}. It follows from Theorem \ref{etathm} that
\begin{equation} \label{etaeq}
    \eta(\gamma)  \;=\;
         \bar{\Phi}(\g,\varphi(\gamma))\cdot\bar{\Psi}(\g,\chi(\gamma)),
\end{equation}
for all $\gamma\in\pi_1\tp(S)$.  Here $\bar{\Phi}(\g,\beta)$ and
$\bar{\Psi}(\g,h)$ are as in \S \ref{locals5}. Therefore, if we know
$\phi$ and $\chi$ explicitly, we can also compute $\eta\tp$.

\begin{rem} \label{practiserem}
  In practice, it is not always so easy to describe an affine frame
  $X\cong\PP^1_S$ and the induced lift $\varphi$ of $\varphib$
  explicitly. In many cases, this is possible, using the methods of
  \cite{Michael}. However, for applications to the Regular Inverse
  Galois Problem, it is usually sufficient to determine the image of
  the projective representation associated to $\eta\tp$, and one can
  proceed as follows.
  
  Let $\varphi:\pi_1\tp(S,s_0)\to A_4$ and $\chi:\pi_1(S,s_0)\to G$ be any
  pair of group homomorphisms such that $\varphi$ is a lift of $\varphib$ and
  such that \eqref{ggeq} holds. (Using the fact that $\pi_1\tp(S,s_0)$ is a
  free group, it is easy to see that such a pair always exists.) The choice of
  $(\varphi,\chi)$ determines a representation $\rho':\pi_1\tp(U)\to\GL_n(R)$
  extending $\rho_0$; it corresponds to a variation $\V'$ of $\V_0\an$. Let
  $\eta':\pi_1\tp(S,s_0)\to\GL(W_\g)$ be the representation corresponding to
  the parabolic cohomology of $\V'$. By Remark \ref{existencerem} (iii), the
  projective representations associated to $\eta\tp$ and $\eta'$ are equal.
  See the next section, in particular \S \ref{example3}.
\end{rem}

%-----------------------------------------------------------------------

\section{An example} \label{example}

We work out one particular example of the construction described in the last
section. In this example, the Hurwitz space is a rational curve. As a result,
we obtain regular realizations over $\QQ(t)$ of certain simple groups
$\PSL_2(\FF_{p^2})$.

\subsection{} \label{example1}

Let $G:=\PSL_2(7)\times\ZZ/3\ZZ$. Given a conjugacy class $C$ of elements of
the group $\PSL_2(7)$, we denote by $C_i$ the conjugacy class of $(g,i)$ in
$G$, where $g\in C$ and $i\in\ZZ/3\ZZ$. The conjugacy classes of $\PSL_2(7)$
are denoted in the standard way (see \cite{Atlas}). For instance, $2a$ is the
unique class of elements of $\PSL_2(7)$ order $2$. Set
\[
     \bC \;:=\; (2a_0,2a_0,3a_1,3a_2).
\]
A computer calculation shows that the set $\Ni\red(\bC)$ has $90$
elements and that the mapping class group $\Bb_4$ acts transitively.
Since $\bC$ is {\em rational} (in the sense of \cite{Voelklein}), the
connected component $S:=H\red(\bC)$ of the Hurwitz space $H_4\red(G)$
corresponding to this orbit is defined over $\QQ$. So $S$ is
a smooth, affine and absolutely irreducible curve over
$\QQ$. Furthermore, our explicit knowledge of the braid action on
$\Ni\red(\bC)$ can be used to show that the complete model $\bar{S}$
of $S$ has genus $0$

\begin{lem} \label{versallem}
\begin{enumerate}
  \item
    The curve $S$ is isomorphic to a dense open subset of $\PP^1_{\QQ}$. 
  \item 
    There exists a versal $G$-Galois cover $f:Y\to X$ over $S$.
\end{enumerate}
\end{lem}

\proof We have to show that $\bar{S}\cong\PP^1_{\QQ}$. Since $\bar{S}$ has
genus $0$, it is well known that it suffices to find a $\QQ$-rational
effective divisor of odd degree on $\bar{S}$. The description of the covering
$j:S_\CC\to\CC$ in terms of the braid action shows that the set of cusps
(i.e.\ the points of $\bar{S}-S$) is such a divisor, of degree $17$. This
finishes the proof of (i).

Let $s=\Spec k\to S$ be a geometric point of $S$ and denote by
$f_s:Y_s\to X_s$ the $G$-cover of type $\bC$ corresponding to $s$. Let
$x_1,\dots,x_4$ denote the branch points of $f_s$, ordered in such a
way that $x_i$ corresponds to the conjugacy class $C_i$. By
definition, we have an injection $G\inj{\rm Aut}_k(Y_s)$.  We claim
that the centralizer of $G$ inside ${\rm Aut}_k(Y_s)$ is equal to the
center of $G$ (which is cyclic of order $3$).  Indeed, suppose that
$\sigma:Y_s\iso Y_s$ is an automorphism which centralizes the action
of $G$.  The automorphism $\sigma':X_s\iso X_s$ induced by $\sigma$
fixes the set $\{x_1,x_2\}$ and the branch points $x_3$ and $x_4$. If
$\sigma'$ were nontrivial, it would be of order $2$, and there would
exist a reduced Nielsen class $[\g]\red\in\Ni\red(\bC)$ which is fixed
by the element $\betab_1\betab_2\betab_1\in\Bb_4$. However, one checks
that such a Nielsen class does not exist, so $\sigma'$ is the
identity. This proves the claim. 

The claim implies that for any $G$-cover $f:Y\to X$ over a scheme $T$ whose
classifying morphism $\varphi_f:T\to H_4(G)$ has its image contained in $S$,
the automorphism group of $f$ is canonically isomorphic to the center of $G$.
It is shown in \cite{diss} that the category of all (families of) $G$-covers
of type $\bC$ is a gerbe over the Hurwitz space $S=H(\bC ).$ In our case, the
band of this gerbe is simply the constant group scheme $\ZZ/3\ZZ.$ By general
results on non-abelian cohomology, the gerbe is represented by a class
$\omega$ in $H^2(S,\ZZ/3\ZZ),$ and the existence of a global section (i.e.,
the neutrality of the gerbe) is equivalent to the vanishing of $\omega.$ See
also \cite{DDE}.

Let $K$ denote the function field of $S$. By (i), $K=\QQ(t)$ is a rational
function field. Since $S$ is affine, we may regard $\omega$ as an element of
the Galois cohomology group $H^2(K,\ZZ/3\ZZ)$. We can give a more concrete
description of $\omega$, as follows. Let $f_{\bar{K}}:Y_{\bar{K}}\to
X_{\bar{K}}$ denote the $G$-cover of type $\bC$ corresponding to the generic
geometric point $\Spec\bar{K}\to S$. For $\sigma\in{\rm Gal}(\bar{K}/K)$, let
$f_{\bar{K}}^\sigma$ denote the conjugate $G$-cover. By definition of the
field $K$, the cover $f_{\bar{K}}^\sigma$ is isomorphic to $f_{\bar{K}}$,
i.e.\ there exists a commutative diagram
\[\begin{CD}
   Y_{\bar{K}} @>{\psi_\sigma}>> Y_{\bar{K}}^\sigma       \\
   @V{f_{\bar{K}}}VV       @V{f_{\bar{K}}^\sigma}VV       \\
   X_{\bar{K}} @>{\varphi_\sigma}>> X_{\bar{K}}^\sigma,    \\
\end{CD}\]
where $\psi_\sigma$ and $\varphi_\sigma$ are $\bar{K}$-linear isomorphisms and
$\psi_\sigma$ is also $G$-equivariant. Note that $\psi_\sigma$ is not uniquely
determined by $\sigma$: we may compose it with an element of the center of
$G$. However, $\varphi_\sigma$ is uniquely determined by $\sigma$ and
therefore satisfies the obvious cocycle relation. We conclude that there
exists a (unique) model $X_K$ of $X_{\bar{K}}$ over $K$ such that
$\varphi_\sigma$ is determined by the isomorphism $X_{\bar{K}}\cong
X_K\otimes\bar{K}$, in the obvious way. In the language of \cite{DebesDouai},
we obtain the following result. The {\em field of moduli} of the $G$-cover
$f_{\bar{K}}$ with respect to the extension $\bar{K}/K$ and the model $X_K$ of
$X_{\bar{K}}$ is equal to $K$. Moreover, the class $\omega\in H^2(K,\ZZ/3\ZZ)$
is the obstruction for $K$ to be a field of definition.
  
The curve $X_K$ is isomorphic to the projective line over $K$ if and only if
it has a $K$-rational point. Moreover, there exists a quadratic extension
$L/K$ such that $X_L:=X_K\otimes L$ has an $L$-rational point and is
isomorphic to $\PP^1_L$. It follows from a theorem of D\`ebes and Douai
\cite{DebesDouai} that $L$ is a field of definition of $f_{\bar{K}}$ (here we
use that the center of $G$ is a direct summand of $G$). In other words, the
restriction of $\omega$ to $L$ vanishes. But by \cite{SerreCG}, Chap. I.2,
Prop. 9, the restriction map $H^2(K,\ZZ/3\ZZ)\to H^2(L,\ZZ/3\ZZ)$ is an
isomorphism. We conclude that $\omega=0$, which finishes the proof of the
proposition.  \Endproof

\begin{rem}
  The lemma shows that there exist infinitely many non-isomorphic
  $G$-covers $f_0:Y_0\to X_0$ defined over $\QQ$. However, we do not
  know whether we can find any such $G$-cover with
  $X_0\cong\PP^1_{\QQ}$. So we do not know whether the lemma produces
  any regular realizations of the group $G$ over $\QQ(t)$. 
\end{rem}

\subsection{} \label{example2}

Let $\g\in\E(\bC)$ be any generating system of type $\bC$; for
instance, we could take
\begin{multline*}
  \g \;:=\; (\;(1,2)(3,4)(5,8)(6,7),\; ( 1, 6)( 2, 5)( 3, 7)( 4, 8), \\
     (1,3,8)(4,5,7)(9,11,10),\; (1,3,7)(2,8,6)(9,10,11)\;)
\end{multline*}
(here we have chosen a faithful permutation representation $G\inj S_{11}$).
The group $G$ admits a faithful and absolutely irreducible linear
representation of dimension $3$, defined over the number field
$K:=\QQ(\sqrt{-3},\sqrt{-7})$. This representation is already defined over
$R:=\OO_K[1/7!]$. From now on, we will consider $G$ as a subgroup of
$\GL_3(R)$. Note that the matrices $g_1,g_2$ are conjugate to the diagonal
matrix ${\rm diag}(1,1,-1)$, and that $g_3$ (resp.\ $g_4$) is conjugate to
${\rm diag}(1,1,\omega)$ (resp.\ ${\rm diag}(1,1,\omega^2)$), with
$\omega:=(-1+\sqrt{-3})/2$. Therefore, Remark
\ref{locals2rem} shows that the $R$-module $W_\g$ defined in \S \ref{locals2}
is locally free of rank $2$.

We can now apply the construction of \S \ref{hurwitz3} to the versal
$G$-cover $f:Y\to X$ of Lemma \ref{versallem}. In particular, for
each prime ideal $\p$ of $K$ with residue characteristic $p\geq 11$,
we obtain a representation
\[
     \eta_\p:\pi_1(S,s_0) \;\To\; \GL_2(\OO_\p).
\]
We let $\eta_\p\geo$ denote the restriction of $\eta_\p$ to the
geometric fundamental group $\pi_1(S_{\Qb})$. We write $\lambda_\p$ and
$\lambda_\p\geo$ for the induced projective representations. We say that
$\eta_\p$ (resp.\ $\lambda_\p$) is {\em regular} if it has the same image
as $\eta_\p\geo$ (resp.\ $\lambda_\p\geo$).

\begin{thm} \label{mainthm}
  Suppose that $p>7$ is not totally split in the extension $K/\QQ$. Then
  $\lambda_\p$ is regular and has image $\PSL_2(\OO_\p)$.
\end{thm}

Before we give the proof of this theorem, let us mention the following
immediate corollary.

\begin{cor} \label{maincor}
  The simple groups $\PSL_2(p^2)$ admit regular realizations over
  $\QQ(t)$, for $p\not\equiv 1,4,16\mod{21}$.
\end{cor}

Note that regular realizations of $\PSL_2(p^2)$ are already known for $p\leq
7$. If $p$ is congruent to $1$, $4$ or $16$ modulo $21$ then our construction
gives a nonregular Galois extension of $\QQ(t)$ with group $\PGL_2(p)$.

\subsection{}  \label{example3}

Set $\Gamma:=\pi_1\tp(S,s_0)$; the first step in the proof of Theorem
\ref{mainthm} is to determine the image of the projective representation
$\lambda\tp:\Gamma\to\PGL(W_\g)$ associated to $\eta\tp:\Gamma\to\GL(W_\g)$.
Since $S$ is isomorphic to the Riemann sphere minus $17$ points, there exist
generators $\gamma_1,\ldots,\gamma_{17}$ of $\Gamma$, subject to the relation
$\prod_j\gamma_j=1$. Our strategy is to explicitly compute the image of
$\gamma_j$ in $\PGL(W_\g)$, for a certain choice of the generators $\gamma_j$.

Let $\varphib:\Gamma\to\Bb_4$ be the group homomorphism induced from
the branch locus configuration $(X,D)$ of the versal $G$-cover $f:Y\to
X$. By construction, there exists a reduced Nielsen class in
$\Ni\red(\bC)$ which is stabilized by the image of $\varphib$. Since
the action of $\Bb_4$ on $\Ni\red(\bC)$ is transitive, we may
normalize things in such a way that the class $[\g]\red$ of our
originally chosen tuple $\g$ is stabilized by $\varphib(\Gamma)$. It
is well known (see e.g.\ \cite{DebesFried}) that there exist generators
$\delta_0,\delta_\infty,\delta_{1728}$ of $\pi_1\tp(\CC-\{0,1728\})$,
with relation $\delta_0\delta_\infty\delta_{1728}=1$,
which are mapped to $\betab_1\betab_2$, $\betab_1$ and
$\betab_1\betab_2\betab_1$, under the natural map
\[
      \pi_1\tp(\CC-\{0,1728\}) \;\To\; \Bb_4.
\]
Let $\Gamma'\subset\pi_1\tp(\CC-\{0,1728\})$ be the inverse image of the
stabilizer of the reduced Nielsen class $[\g]\red$. We may identify $\Gamma'$
with the fundamental group of $S':=j^{-1}(\CC-\{0,1728\})\subset S$. It is a
straightforward, although combinatorially involved problem to write down a
list of generators of the free group $\Gamma'$, given as words in the
generators $\delta$. Moreover, one can choose these generators in such a way
that the usual product-$1$-relation holds and that each of them represents a
simple closed loop around one of the points missing from $S'$. Let
$\gamma_1,\ldots,\gamma_{17}\in\Gamma$ be those generators representing a loop
around a cusp (i.e.\ a point $s\in\bar{S}$ with $j(s)=\infty$). Note that
$\gamma_j$ is conjugate (inside the group $\pi_1\tp(\CC-\{0,1728\})$) to a
certain power of $\delta_\infty$. The other generators, representing a loop
around one of the points of $S-S'$, are conjugate either to $\delta_0^3$ or to
$\delta_{1728}^2$, so their image in $\Bb_4$ is $1$. It follows that the map
$\Gamma'\to\Bb_4$ factors over the natural, surjective map $\Gamma'\to\Gamma$.
Denoting the image of $\gamma_j$ in $\Gamma$ by the same name, we have found
explicit generators $\gamma_1,\ldots,\gamma_{17}$ of $\Gamma$, with relation
$\prod_j\gamma_j=1$, and their images under the map $\varphib:\Gamma\to\Bb_4$.

It is easy to find, for all $j=1,\ldots,17$, an element $\gamma_j'\in
A_4$ which lifts $\varphib(\gamma_j)$ and an element $h_j\in G$ such
that
\[
           \g^{\gamma_j'} \;=\; \g^{h_j}.
\]
Moreover, we may do this in such a way that $\prod_j\gamma_j'=1$ and
$\prod_jh_j=1$. In other words, we can choose homomorphisms $\varphi:\Gamma\to
A_4$ and $\chi:\Gamma\to G$ as in Remark \ref{practiserem}. In fact, the lift
$\varphi$ is unique, because the Klein four group $Q$ acts faithfully on
$\Ni(\CC)$. On the other hand, $\chi$ is only determined up to multiplication
of $h_j=\chi(\gamma_j)$ by a central element of order $3$. This corresponds to
the fact that the versal $G$-cover $f:Y\to X$ over $S$ may be twisted by
characters of order $3$. (It is not clear how to find $\chi$ corresponding to
a versal cover $f$ defined over $\QQ$). By formula \eqref{etaeq} and Remark
\ref{practiserem} we have
\[
     \eta\tp(\gamma_j) \;=\; c_j\cdot
          \bar{\Phi}(\g,\gamma_j')\cdot\bar{\Psi}(\g,h_j),
\]
for some skalar $c_j\in K^\times$. (In fact, $c_j$ is a third root of unity and
we have $\prod_jc_j=1$.) Set
$b_j:=\bar{\Phi}(\g,\gamma_j')\cdot\bar{\Psi}(\g,h_j)$. By construction, $b_j$
is an invertible $2$-by-$2$-matrix with entries in $R$ such that $\prod_j
b_j=1$.

Using a computer program written in {\em GAP}, the authors have computed the
matrices $b_j$ explicitly. It turns out that $12$ of the $b_j$ are
transvections and $5$ are homologies with eigenvalues $1,\omega$ or
$1,\omega^2,$ where $\omega$ denotes a primitive third root of unity (see
\cite{DettwReiterKatz} for notations). One finds that the trace of the matrix
$b_1b_2$ is a generator of the extension $K/\QQ$. Moreover, one checks that
for every prime $p>7$ one can find a pair of transvections $b_i$, $b_j$ whose
commutator is not congruent to the identity, modulo any prime ideal $\p$ above
$p$. This information suffices to show that for a prime $p>7$ which is not
totally split in $K/\QQ$, the image of the residual projectivized
representation $\bar{\lambda}_\p\geo$ associated to $\eta_\p\geo$ is equal to
$\PSL_2(p^2)$ (we may identify the residue field of $\p$ with $\FF_{p^2}$).
By a well known argument (see e.g.\cite{SerreAbelian}), it follows that the
image of the projective representation $\lambda_\p\geo$ associated to
$\eta_\p\geo$ is equal to $\PSL_2(\OO_\p)$.

The only thing left to prove is that the image of the full projective
representation $\lambda_\p$ is equal to $\PSL_2(\OO_\p)$ as well. Again, it
suffices to show that the image of the residual projectivized representation
$\bar{\lambda}_\p$ is equal to $\PSL_2(p^2)$. 

One observes that there are exactly five ramification points
$s_1,\ldots,s_5\in\bar{S}$ of the map $j:\bar{S}\to\PP^1$ above $\infty$ whose
ramification index is equal to $4$. One also observes that the matrices
$b_{j_\mu}$ corresponding to the points $s_\mu$ are transvection for
$\mu=1,\ldots,4$, whereas $b_{j_5}$ is a homology. It follows that the set
$\{s_1,\ldots,s_4\}$ is rational, i.e.\ fixed by the action of
$\Gal(\Qb/\QQ)$. Furthermore, the transvections $b_{\mu_1},\ldots,b_{\mu_4}$
are all conjugate to each other by elements of $\SL_2(K)$. One concludes that
the image of these transvections give rise to conjugate transvections in the
image of $\etab_\p\geo$. The conjugacy class of these transvections is a
rational class, in the sense of \cite{Voelklein}. 

Suppose that there exists an element $\sigma\in\Gal(\Qb/\QQ)$ and a lift
$\alpha\in\pi_1(S)$ of $\sigma$ such that $\etab_\p(\alpha)\in\GL_2(p^2)$ does
not lie in $\SL_2(p^2)Z(\GL_2(p^2))$. Using the branch cycle argument (as in
the proof of \cite{VoelkleinCrelle}, Corollary 4.6), one would conclude that
$\sigma$ does not fix the set $\{s_1,\ldots,s_4\}$. But this would be a
contradiction to the assertion made above. It follows that the image of
$\bar{\lambda}_\p$ is equal to $\PSL_2(p^2)$. The proof of Theorem
\ref{mainthm} is now complete.  \Endproof

\begin{rem}
  It is also possible to prove Theorem \ref{mainthm} by using a
  generalisation of methods developed by V\"olklein, see e.g.\ 
  \cite{VoelkleinCrelle}. Here is a brief outline. One constructs a
  certain projective system of finite \'etale covers of $S$ (these
  covers are themselves Hurwitz spaces). Our representation $\eta\geo$
  essentially corresponds to the inverse limit of the Galois closures
  of these covers. One then has to show, using the theory of
  \cite{FriedVoe91}, that this projective system has $\QQ$ as a field
  of moduli. One deduces that $\eta\geo$ has $\QQ$ as a field of
  linear moduli (in our approach this is automatic, as $\eta\geo$ is
  the restriction of a representation $\eta$ of the full arithmetic
  fundamental group). The rest of the proof goes as above.
\end{rem}

%-------------------------------------------------------------------------

\section{The monodromy of the Picard--Euler system} \label{picard}

\subsection{} \label{picard1}

In this section, we use the notation of \S \ref{hurwitz}. However, all
varieties are defined over the complex numbers. Let $G$ be a cyclic
group of order $3$, with generator $\sigma$. We fix a nontrivial
character $\chi:G\inj\CC^\times$. Then $\omega:=\chi(\sigma)$ is a
primitive third root of unity. We consider the generating system
\[
  \g \;=\; (\sigma,\sigma,\sigma,\sigma,\sigma^2) \;\in\;\E_5(G).
\]
The orbit $O$ of $\g$ under the action of the braid group consists simply
of the five permutations of $\g$. Let 
\[
    S \;:=\; \{\,(s,t)\in\CC^2 \;\mid\; s,t\not=0,1,\; s\not=t\;\},
\]
and let $X:=\PP^1_S$ denote the relative projective line over $S$. 
The equation 
\begin{equation} \label{picardeq}
       y^3 \;=\; x(x-1)(x-s)(x-t)
\end{equation}
defines a finite Galois cover $f:Y\to X$ of smooth projective curves
over $S$, tamely ramified along the divisor
$D:=\{0,1,s,t,\infty\}\subset X$. We identify the Galois group of $f$
with $G$ in such a way that $\sigma^*y=\omega\cdot y$. The classifying
map $\varphi_f:S\to H_5\red(G)$ is a finite \'etale cover of the
connected component $H(O)\subset H_5\red(G)$ corresponding to the
braid orbit $O$ of $\g$. Thus, in the terminology of \S \ref{hurwitz},
$f:Y\to X$ is a versal family of $G$-covers of type $\g$. 
A $G$-cover of $\PP^1$ of type $\g$ is called a {\em
  Picard curve}, see e.g.\ \cite{HolzapfelBall}.  

Let $K:=\QQ(\omega)$ denote the field of third roots of unity and
$\OO_K=\ZZ[\omega]$ its ring of integers. The family of $G$-covers
$f:Y\to X$ together with the character $\chi$ of $G$ give rise to a
local system of $\OO_K$-modules on $U:=X-D$, see \S \ref{hurwitz3}. Set
$s_0:=(2,3)\in S$ and let $\V_0$ denote the restriction of $\V$ to the
fibre $U_0=\AA^1_\CC-\{0,1,2,3\}$ of $U\to S$ over $s_0$.  We consider
$\V$ as a variation of $\V_0$ over $S$. Let $\W$ denote the parabolic
cohomology of this variation; it is a local system of $\OO_K$-modules
of rank three, see Remark \ref{locals2rem}.  Let
$\chi':G\inj\CC^\times$ denote the conjugate character to $\chi$ and
$\W'$ the parabolic cohomology of the variation of local systems $\V'$
corresponding to the $G$-cover $f$ and the character $\chi'$. We write
$\W_\CC$ for the local system of $\CC$-vectorspaces $\W\otimes\CC$. The maps
$\pi_Y:Y\to S$ and $\pi_X:X\to S$ denote the natural projections.

\begin{prop} \label{picardprop1}
  We have a canonical isomorphism of local systems
  \[
       R^1\sing\pi_{Y,*}\underline{\CC} \;\cong\; 
            \W_\CC \,\oplus\, \W'_\CC.
  \]
  This isomorphism identifies the fibres of $\W_\CC$ with the
  $\chi$-eigenspace of the singular cohomology of the Picard curves
  of the family $f$.  
\end{prop}

\proof
The group $G$ has a natural left action on the sheaf $f_*\underline{\CC}$.
It is easy to that we have a canonical isomorphism of sheaves on $X$
\[
       f_*\underline{\CC} \;\cong\; \underline{\CC} \,\oplus\,
        j_*\V_\CC \,\oplus\, j_*\V',
\]
which identifies $j_*\V_\CC$, fibre by fibre, with the
$\chi$-eigenspace of $f_*\underline{\CC}$. Now the Leray spectral
sequence for the composition $\pi_Y=\pi_X\circ f$ gives
isomorphisms of sheaves on $S$
\[
   R^1\pi_{Y,*}\underline{\CC} \;\cong\;
     R^1\pi_{X,*}(f_*\underline{\CC}) \;\cong\;
       \W_\CC \,\oplus\, \W'_\CC.
\]
Note that $R^1\pi_{X,*}\underline{\CC}=0$ because the genus of $X$ is zero.
Since the formation of $R^1\pi_{Y,*}$ commutes with the $G$-action, the
proposition follows.  \Endproof

% \begin{rem} \label{picardrem1}
% \begin{enumerate}
% \item The $G$-action on $Y$ induces an $\OO_K$-multiplication on the
%   relative Jacobian $J_{Y/S}$. The associatation $Y\;\mapsto\;
%   J_{Y/S}$ gives rise to a map from the Hurwitz space $H_5\red(G)$
%   (i.e.\ the moduli space of Picard curves) to a certain Shimura
%   variety of type $\mathbb{U}(2,1)$, see \cite{HolzapfelBall}.
% \item Remark (i) above generalizes to the more general situation of an
%   arbitry finite group $G$ with a faithful representation
%   $G\inj\GL_n(K)$, see \cite{LangeRecillas}. This gives a quite
%   general constructions of maps from Hurwitz spaces to Shimura
%   varieties of PEL-type, which we intend to study in a future paper. 
% \end{enumerate}
% \end{rem}

\subsection{}   \label{picard2}

The comparison theorem between singular and deRham cohomology
identifies $R^1\sing\pi_*\underline{\CC}$ with the local system of
horizontal sections of the relative deRham cohomology module
$R^1\dR\pi_*\OO_Y$, with respect to the Gauss-Manin connection. The
$\chi$-eigenspace of $R^1\dR\pi_*\OO_Y$ gives rise to a Fuchsian
system known as the Picard--Euler system. In more classical terms, the
Picard--Euler system is a set of three explicit partial differential
equations in $s$ and $t$ of which the period integrals
\[
      I(s,t;a,b) \;:=\; 
        \int_a^b\,\frac{{\rm d}\,x}{\sqrt[3]{x(x-1)(x-s)(x-t)}}
\]
(with $a,b\in\{0,1,s,t,\infty\}$) are a solution. See \cite{Picard83},
\cite{HolzapfelEuler}, \cite{HolzapfelBall}. It follows from
Proposition \ref{picardprop1} that the monodromy of the Picard--Euler
system can be identified with the representation
$\eta:\pi_1(S)\to\GL_3(\OO_K)$ corresponding to the local system $\W$. 

\begin{thm}[Picard] \label{picardthm}
  For suitable generators $\gamma_1,\ldots,\gamma_5$ of the fundamental
  group $\pi_1(S)$, the matrices
  $\eta(\gamma_1),\ldots,\eta(\gamma_5)$ are equal to 
  \begin{gather*}
     \begin{pmatrix}
        \omega^2        & \;0\;          & 1-\omega            \\
        \omega-\omega^2 & 1              & \omega^2 -1         \\
        0               & 0              & 1 
     \end{pmatrix},\;
     \begin{pmatrix}
        \omega^2        & \;0\;          & 1-\omega^2          \\
        1-\omega^2      & 1              & \omega^2-1          \\
        0               & 0              & 1 
     \end{pmatrix},\;
     \begin{pmatrix}
        \;\;1\;\;       & 0              & 0                    \\
        0               & \omega         & \omega^2-1           \\
        0               & \omega^2-1     & -2\omega   
     \end{pmatrix},\\
     \begin{pmatrix}
        \;\;\omega^2\;  & \;\;0\;\;      & \;\;0\;\;            \\
        0               & 1              & 0                    \\
        0               & 0              & 1 
     \end{pmatrix},\;
     \begin{pmatrix}
        \omega^2        & \omega-\omega^2& \;0\;\\
        0               & 1              & 0    \\
        1-\omega        & \omega^2-1     & 1 
     \end{pmatrix}.
  \end{gather*}
\end{thm}     

\begin{figure}
\begin{center}

\setlength{\unitlength}{0.0011in}
\begingroup\makeatletter\ifx\SetFigFont\undefined%
\gdef\SetFigFont#1#2#3#4#5{%
  \reset@font\fontsize{#1}{#2pt}%
  \fontfamily{#3}\fontseries{#4}\fontshape{#5}%
  \selectfont}%
\fi\endgroup%
{\renewcommand{\dashlinestretch}{30}
\begin{picture}(4224,899)(0,200)
\path(12,864)(612,864)
\path(12,114)(612,114)
\path(87,864)(87,114)
\path(237,864)(237,114)
\drawline(537,114)(537,114)
\path(387,114)(387,118)(387,126)
        (387,139)(387,159)(387,183)
        (387,211)(387,240)(387,269)
        (387,296)(387,322)(387,345)
        (387,365)(387,384)(387,400)
        (387,414)(387,427)(387,439)
        (387,458)(387,475)(387,491)
        (387,507)(387,523)(387,538)
        (387,551)(387,559)(387,563)(387,564)
\path(541,110)(541,111)(540,115)
        (538,125)(534,141)(530,160)
        (525,180)(520,199)(515,216)
        (511,230)(507,243)(503,254)
        (498,265)(493,277)(486,289)
        (478,302)(468,317)(457,333)
        (445,349)(434,364)(425,374)
        (422,379)(421,380)
\path(353,440)(352,441)(350,445)
        (344,453)(338,464)(332,473)
        (328,482)(325,489)(322,496)
        (321,502)(320,508)(319,515)
        (319,522)(320,529)(322,537)
        (325,545)(329,552)(334,560)
        (340,567)(347,575)(355,583)
        (365,591)(377,600)(389,610)
        (401,619)(414,629)(425,638)
        (437,648)(447,656)(457,666)
        (466,675)(475,685)(484,696)
        (492,707)(500,718)(507,729)
        (512,740)(517,750)(521,760)
        (524,770)(527,780)(529,792)
        (531,805)(533,821)(535,837)
        (536,852)(537,862)(537,867)(537,868)
\path(387,673)(387,868)
\path(987,860)(987,110)
\path(912,864)(1512,864)
\path(912,114)(1512,114)
\path(1441,114)(1441,115)(1438,118)
        (1433,126)(1424,139)(1414,154)
        (1402,170)(1391,185)(1381,199)
        (1371,210)(1361,221)(1352,230)
        (1343,239)(1332,247)(1321,256)
        (1308,265)(1293,275)(1276,285)
        (1256,297)(1235,309)(1215,321)
        (1198,330)(1188,336)(1183,339)(1182,339)
\path(1812,114)(2412,114)
\path(1812,864)(2412,864)
\drawline(2337,114)(2337,114)
\path(2337,118)(2337,119)(2335,123)
        (2332,131)(2327,144)(2322,157)
        (2315,170)(2309,182)(2303,193)
        (2295,203)(2287,213)(2279,222)
        (2270,232)(2260,242)(2249,253)
        (2236,264)(2222,275)(2208,287)
        (2194,298)(2179,309)(2165,319)
        (2151,328)(2137,337)(2124,344)
        (2111,351)(2098,358)(2083,366)
        (2066,374)(2047,382)(2027,392)
        (2006,401)(1985,410)(1967,418)
        (1953,424)(1944,427)(1941,429)(1940,429)
\path(1853,463)(1852,464)(1847,469)
        (1839,477)(1832,485)(1828,492)
        (1825,499)(1823,505)(1822,512)
        (1822,519)(1822,527)(1824,535)
        (1827,544)(1831,551)(1836,559)
        (1841,565)(1848,571)(1856,577)
        (1866,583)(1878,590)(1892,597)
        (1907,603)(1924,609)(1941,615)
        (1961,621)(1974,625)(1988,629)
        (2003,633)(2019,637)(2036,641)
        (2054,646)(2072,650)(2091,655)
        (2109,660)(2127,664)(2145,669)
        (2162,674)(2177,678)(2192,683)
        (2205,688)(2218,692)(2236,700)
        (2252,708)(2267,716)(2280,726)
        (2292,735)(2301,745)(2309,754)
        (2316,763)(2321,772)(2325,781)
        (2328,791)(2331,803)(2333,816)
        (2335,832)(2336,849)(2337,863)
        (2337,871)(2337,872)
\path(1887,557)(1887,114)
\path(1887,632)(1887,864)
\path(2041,857)(2041,673)
\path(2037,602)(2037,429)
\path(2033,347)(2033,118)
\path(2187,857)(2187,718)
\path(2191,650)(2191,343)
\path(2187,260)(2187,118)
\path(2712,864)(3312,864)
\path(2712,114)(3312,114)
\path(2787,864)(2787,114)
\path(3237,864)(3237,114)
\path(3091,118)(3091,119)(3090,123)
        (3088,133)(3084,149)(3080,168)
        (3075,188)(3070,207)(3065,224)
        (3061,238)(3057,251)(3053,262)
        (3048,273)(3043,285)(3036,297)
        (3028,310)(3018,325)(3007,341)
        (2995,357)(2984,372)(2975,382)
        (2972,387)(2971,388)
\path(2937,665)(2937,860)
\path(2903,440)(2902,441)(2900,445)
        (2894,453)(2888,464)(2882,473)
        (2878,482)(2875,489)(2872,496)
        (2871,502)(2870,508)(2869,515)
        (2869,522)(2870,529)(2872,537)
        (2875,545)(2879,552)(2884,560)
        (2890,567)(2897,575)(2905,583)
        (2915,591)(2927,600)(2939,610)
        (2951,619)(2964,629)(2975,638)
        (2987,648)(2997,656)(3007,666)
        (3016,675)(3025,685)(3034,696)
        (3042,707)(3050,718)(3057,729)
        (3062,740)(3067,750)(3071,760)
        (3074,770)(3077,780)(3079,792)
        (3081,805)(3083,821)(3085,837)
        (3086,852)(3087,862)(3087,867)(3087,868)
\path(2937,114)(2937,118)(2937,126)
        (2937,139)(2937,159)(2937,183)
        (2937,211)(2937,240)(2937,269)
        (2937,296)(2937,322)(2937,345)
        (2937,365)(2937,384)(2937,400)
        (2937,414)(2937,427)(2937,439)
        (2937,458)(2937,475)(2937,491)
        (2937,507)(2937,523)(2937,538)
        (2937,551)(2937,559)(2937,563)(2937,564)
\path(1103,384)(1102,385)(1097,390)
        (1090,399)(1083,408)(1079,416)
        (1076,425)(1074,433)(1073,442)
        (1072,452)(1072,463)(1072,474)
        (1074,485)(1076,494)(1078,502)
        (1081,510)(1085,518)(1091,525)
        (1098,532)(1106,539)(1115,546)
        (1125,553)(1136,559)(1144,564)
        (1153,569)(1163,575)(1173,581)
        (1185,587)(1197,594)(1209,601)
        (1222,608)(1234,615)(1246,622)
        (1258,629)(1269,636)(1280,643)
        (1291,651)(1302,659)(1314,667)
        (1325,676)(1336,685)(1347,695)
        (1357,704)(1366,713)(1375,723)
        (1382,732)(1388,741)(1395,752)
        (1401,764)(1406,777)(1411,793)
        (1416,811)(1421,831)(1425,848)
        (1428,861)(1429,867)(1429,868)
\path(1137,857)(1137,598)
\path(1137,519)(1137,114)
\path(1287,864)(1287,684)
\path(1287,602)(1287,324)
\path(1287,241)(1287,118)
\drawline(3612,864)(3612,864)
\path(3612,864)(4212,864)
\path(3612,114)(4212,114)
\path(4137,864)(4137,114)
\path(3660,380)(3659,381)(3654,386)
        (3647,395)(3640,404)(3636,412)
        (3633,421)(3631,429)(3630,438)
        (3629,448)(3629,459)(3629,470)
        (3631,481)(3633,490)(3635,498)
        (3638,506)(3642,514)(3648,521)
        (3655,528)(3663,535)(3672,542)
        (3682,549)(3693,555)(3701,560)
        (3710,565)(3720,571)(3730,577)
        (3742,583)(3754,590)(3766,597)
        (3779,604)(3791,611)(3803,618)
        (3815,625)(3826,632)(3837,639)
        (3848,647)(3859,655)(3871,663)
        (3882,672)(3893,681)(3904,691)
        (3914,700)(3923,709)(3932,719)
        (3939,728)(3945,737)(3952,748)
        (3958,760)(3963,773)(3968,789)
        (3973,807)(3978,827)(3982,844)
        (3985,857)(3986,863)(3986,864)
\path(3991,114)(3991,115)(3988,118)
        (3983,126)(3974,139)(3964,154)
        (3952,170)(3941,185)(3931,199)
        (3921,210)(3911,221)(3902,230)
        (3893,239)(3882,247)(3871,256)
        (3858,265)(3843,275)(3826,285)
        (3806,297)(3785,309)(3765,321)
        (3748,330)(3738,336)(3733,339)(3732,339)
\path(3687,523)(3687,118)
\path(3687,857)(3687,598)
\path(3837,856)(3837,676)
\path(3837,591)(3837,313)
\path(3837,241)(3837,118)
\put(50,-30){$0$}
\put(206,-30){$1$}
\put(515,-30){$t$}
\put(365,-30){$s$}
\end{picture}
}

\end{center}
\caption{\label{zopfbild} The braids $\gamma_1,\ldots,\gamma_5$}
\end{figure}
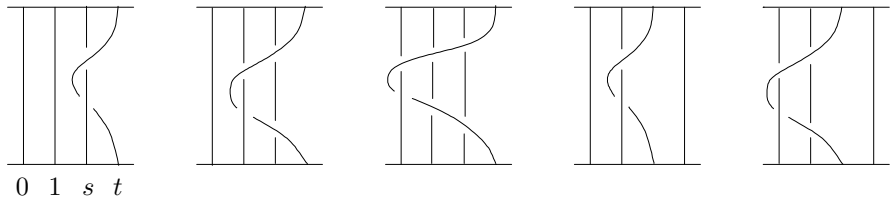

\proof Recall that $A_r$ (resp.\ $B_r$) denotes the Artin (resp.\ 
Hurwitz) braid group on $r$ strands. We identify $A_{r-1}$ with the
subgroup of $B_r$ generated by the standard braids
$\beta_1,\ldots,\beta_{r-2}$. By the results of \S \ref{variation} and
\S \ref{hurwitz}, the representation
$\eta:\pi_1(S)\to\GL(W_\g)\cong\GL_3(\OO_K)$ factors through the map
$\varphi:\pi_1(S)\to A_4\subset B_5$ induced by the branch divisor
$D\subset\PP^1_S$ (we write the points of $D$ in the order
$(0,1,s,t,\infty)$). Using standard methods (see e.g.\
\cite{Voelklein01} or \cite{DettwReiterKatz}), one can show that the
image of $\varphi$ is indeed generated by the five braids
\[
   \beta_3^2,\;\; \beta_3\beta_2^2\beta_3^{-1},\;\;
   \beta_3\beta_2\beta_1^2\beta_2^{-1}\beta_3^{-1},\;\;
   \beta_2^2,\;\; \beta_2\beta_1^2\beta_2^{-1},
\]
see Figure \ref{zopfbild}.  It is clear that these five
braids can be realized as the image under the map $\varphi$ of
generators $\gamma_1,\ldots,\gamma_5\in\pi_1(S)$.

Let $\rho:\pi_1(U)\to G\subset K^\times$ denote the representation
corresponding to the $G$-cover $f:Y\to X$, and $\rho_0:\pi_1(U_0)\to
G$ its restriction to the fibre above $s_0$. Considering the
$\infty$-section as a `tangential base point' for the fibration $U\to
S$, we obtain a section $\pi_1(S)\to\pi_1(U)$. We use this section to
identify $\pi_1(S)$ with a subgroup of $\pi_1(U)$. Let
$\alpha_1,\ldots,\alpha_5$ be the standard generators of $\pi_1(U_0)$.
Using \eqref{picardeq} one checks that $\rho_0$ corresponds to the
tuple $\g=(\sigma,\sigma,\sigma,\sigma,\sigma^2)$, i.e.\ that
$\rho_0(\alpha_i)=g_i$. Also, since the leading coefficient of the
right hand side of \eqref{picardeq} is one, the restriction of $\rho$
to $\pi_1(S)$ is trivial. Hence, by Theorem \ref{etathm}, we have
\[
       \eta(\gamma_i) \;=\; \bar{\Phi}(\g,\varphi(\gamma_i)).
\]
A straightforward computation, using \eqref{locals5eq6} and the
cocycle rule \eqref{locals5eq7}, gives the value of $\eta(\gamma_i)$
(in form of a three-by-three matrix depending on the choice of a basis
of $W_\g$). For this computation, it is convenient to take the classes
of $(1,0,0,0,-\omega^2)$, $(0,1,0,0,-\omega)$ and
$(0,0,1,0,-1)$ as a basis. In order to obtain the matrices
stated in the theorem, one has to use a different basis, i.e.\ conjugate
with the matrix
\[
     B \;=\; \begin{pmatrix} 
               0          & -\omega-1  & -\omega    \\
               \omega+1   & \omega+1   & \omega+1   \\
               1          & 0          & 0          \\
             \end{pmatrix}.
\]
\Endproof

\begin{rem}
\begin{enumerate}
\item Theorem \ref{picardthm} is due to Picard, see \cite{Picard83}, p.\ 125,
  and \cite{Picard84}, p.\ 181. He obtains exactly the matrices given above,
  but he does not list all of the corresponding braids. Holzapfel in
  \cite{HolzapfelEuler} gives a list of five braids which generated
  $\pi_1(S)$, see \cite{HolzapfelEuler}, p.\ 125. But contrary to what is
  claimed, these braids do not correspond to the five matrices found by
  Picard.
\item Essentially the same computation as in the proof of Theorem
  \ref{picardthm}, but for arbitrary cyclic covers of $\PP^1$ and with
  $p$-adic coefficients, can be found in \cite{Voelklein95}. 
\end{enumerate}
\end{rem}

%-------------------------------------------------------------------------

%\bibliographystyle{plain} \bibliography{dw}

\vspace{4ex}

\begin{minipage}[t]{7cm}
IWR, Universit\"at Heidelberg\\
INF 368\\
69120 Heidelberg\\
michael.dettweiler@iwr.uni-heidelberg.de\\
\end{minipage}
\hfill
\begin{minipage}[t]{4.5cm}
\begin{flushright}
Mathematisches Institut\\ 
Universit\"at Bonn\\
Beringstr. 1, 53115 Bonn\\
wewers@math.uni-bonn.de
\end{flushright}
\end{minipage}

\end{document}